\documentclass[11pt]{article}
\usepackage{amsmath, amscd,amsbsy, amsthm, amsfonts, amssymb}
\usepackage{rotating}
\usepackage{dcolumn,multicol,url}

 \theoremstyle{definition}

 \numberwithin{equation}{section}

 \newcommand{\A}{\mathcal{A}}

\def\C{{\mathbf C}}
\def\D{{\mathbf D}}

\def\R{{\mathbb R}}

\def\Z{{\mathbb Z}}
\def\ZZ{{\mathbb Z}}

\def\RR{{\mathbb R}}

\def\Zee{\mathbb{Z}}

\def\Ar{\mathbb{R}}

\def\ZZ{{\mathbb Z}}

\def\RR{{\mathbb R}}

\def\be{\begin{equation}}
\def\ee{\end{equation}}
\def\ba{\begin{eqneqnarray}}
\def\ea{\end{eqneqnarray}}
\def\tilde{\widetilde}

\def\e1{\epsilon}
\def\AAl{\mathcal{A}_{\lambda}}
\def\A0{\stackrel{\circ}{\AAl}}

\def\o1{\widetilde{\Omega}}
\def\01{\widetilde{\Omega}}
\def\c1{\gamma}

\def\g1{\Sigma}

\def\l1{\Lambda}

\def\v1{\varphi}

\def\d1{\delta}

\def\f1{\frac}
\def\t1{\theta}
\def\b1{\beta}
\def\bar{\overline}
\def\bs{\begin{eqneqnarray*}}
\def\es{\end{eqneqnarray*}}
\def\m1{\Theta}
\def\w1{\wedge}

\begin{document}

\title{The Pontrjagin Dual of 4-Dimensional Spin Bordism}

\author{Greg Brumfiel and John Morgan}

\maketitle

\tableofcontents

\newpage

\section{Introduction}

\subsection{Preview}

The goal of this paper is to study the Pontrjagin dual of (reduced) $4$-dimensional Spin bordism. That is to say, we consider the functor from the category of topological spaces to the category of compact abelian groups that associates to each space $X$ the compact group ${\rm Hom}(\widetilde{\Omega}_4^{\rm spin}(X),\Ar/\Zee)$.\footnote{\label{circle} We identify $\R/ \Z$  with the  unit circle $S^1$ via complex exponentiation $e^{2\pi it}$.  For our purposes,  $\R/ \Z$ is advantageous for the values of characters as it allows for simultaneous additive notation for cochain groups with values in $\Z, \Z/n$, and $\R / \Z$.} By $\widetilde{\Omega}_4^{spin}(X)$ we will mean the  reduced bordism, which means bordism classes of maps $M^4 \to X$ where $M^4$ is a Spin manifold which is a Spin boundary.  Equivalently  $M^4$ is a Spin manifold with signature 0. In a previous paper, see [1], we studied the analogous problem for $3$-dimensional Spin bordism. Our work was motivated by some questions from physics, see especially  [4] and [5].  The physicists are primarily interested in the case $X = B\Gamma$, the classifying space of a finite group. We discuss some aspects of the physicists interest in \S1.2. But with our methods there is no reason to specialize to $B\Gamma$, the discussion applies equally  well to the Pontrjagin dual of the Spin bordism of any space $X$.\\

Standard spectral sequence methods of algebraic topology can be used to compute $n$-dimensional Spin bordism of many spaces $X$ in a large range of dimensions $n$.  The answers are typically presented as direct sums of cyclic abelian groups. See for example [3] for computations essentially equivalent to computing Spin bordism groups in many dimensions for many $B\Gamma$.  But physicists sometimes want a more concrete description, produced by local cellular data. We will explicitly construct a group $G^4(X)$ whose elements are equivalence classes of certain triples $(w, p, a)$ of cochains on $X$.  Given a representative of a Spin bordism class $f\colon M^4 \to X$ and a triangulation of $M$ so that $f$ is simplicial, we then produce from $f$ and a suitable  triple of cochains    on $X$ an invariant in $\R/ \Z$ that explicitly realizes a dual pairing $ \widetilde{\Omega}_4^{spin}(X) \times G^4(X) \to \R / \Z$.\\

While, in principle, one could study the Pontrjagin dual of $n$-dimensional Spin bordism for any $n$ by our methods, it turns out that for all $n>4$ the cochain computations required by our methods become too complicated to carry out in practice. Indeed, one might level the same complaint about the $4$-dimensional case we consider here. But because dimensions $3$ and $4$ are  important dimensions for physical applications, and because the $4$-dimensional case is just on the edge of what one can do  by our methods,  we felt it was worth the effort.  Our work is not meant to be a useful computational answer, but rather a concrete description of how one can express a homotopy invariant like the Pontrjagin dual of reduced 4-dimensional Spin bordism directly in terms of a triangulation of $X$.  It turns out this is the same as describing homotopy classes of maps from $X$ to a certain three stage Postnikov tower in direct simplicial terms.  \\

Some of the proofs in this paper follow the lines of proofs  in our 3-dimensional work given in [1]. This is especially true of \S6 and \S7, where the explicit cochain construction of the group $G^4(X)$ is carried out. Also, the essential details of at least the definition of the pairing $ \widetilde{\Omega}_4^{spin}(X) \times G^4(X) \to \R / \Z$ in \S3  follows the lines  of our 3-dimensional work.  But there are some significant differences in the two cases.  In the 3-dimensional work we also began with certain triples of cochains, and we defined a product operation and an equivalence relation on these triples.  We proved directly that the equivalence classes formed an abelian group $G^3(X)$.  Then we defined an explicit map  $G^3(X) \to {\rm Hom}(\Omega^{spin}_3(X), \R / \Z)$.\footnote{$\ \Omega^{spin}_3(pt) = 0$, so reduced and unreduced Spin bordism agree in dimension three.}   It was quite challenging to prove this map was additive  in the $G^3(X)$ variable.  We then proved the map to be a group isomorphism by comparing natural filtrations on both sides.  A simple Postnikov tower existed behind the scenes but we did not use that in any essential way.  \\

In our 4-dimensional work, we emphasize from the outset the known Postnikov tower that represents the Pontrjagin dual of Spin bordism.  Thus, we know homotopy classes of maps from a space $X$ to this tower form the group we want to study. Also, we emphasize that the Spin bordism groups $\widetilde{\Omega}_n^{spin}(X)$ collectively, that is, for all $n$, are part of a general homology theory, so the representing Postnikov towers of the Pontrjagin dual cohomology theory arise from an $\Omega$-spectrum.  In particular,  there is a product on a specific choice of Postnikov tower representing $G^4(X)$ that is homotopy equivalent to a loop space product.  What we then  do is find cochain formulas for the product of triples and the relations among triples that describe the group of homotopy classes of maps to this tower.  There is also essentially only one possible identification of this group with the Pontrjagin dual of Spin bordism that satisfies certain basic constraints. Thus, once we find suitable formulas, the Postnikov tower theory guarantees our pairing $ \widetilde{\Omega}_4^{spin}(X) \times G^4(X) \to \R / \Z$ is a bi-additive dual pairing. Finding some of these cochain formulas in the 4-dimensional case required extensive computer calculations.\\

Although certainly rather complicated, the 3-dimensional treatment in [1] is considerably more leisurely than the 4-dimensional treatment in this paper.  Because there are many parallel steps in the two papers, it might be advisable for anyone looking at this paper to also look at [1].\\

We thank Yi Sun of Columbia University for describing to us how to do the relevant computer computations with Sage and Python.

\subsection{Some Specific Comments Related to Physics}

The physicists are not interested in the dual to Spin bordism of all spaces $X$.  Their specific interest is the case $X = B\Gamma$, the classifying space of a finite group. A map $M \to  B\Gamma$ corresponds to a principal $\Gamma$-bundle covering space of $M$. Our cochain data on $B\Gamma$ translates  to some kind of local information on $M$ with physical significance.  The physicists then want numerical invariants of some kind for this data, which have some bordism invariance properties and a physical interpretation. Our pairing $ \widetilde{\Omega}_4^{spin}(B\Gamma) \times G^4(B\Gamma) \to \R / \Z$ produces all such invariants.  \\

In the relevant problems of physics, one is considering discrete systems, for example, on a lattice in ordinary space or on a triangulation of a manifold. The systems in question have a gauge field with discrete gauge group, i.e., a principal $\Gamma$-bundle for some finite group $\Gamma$, and fermion fields. As a consequence, the theories only make sense when the generalized space-time has a Spin structure. The physical question is to describe invariants for the limiting topological theory as the scale size of the discrete structure (lattice or triangulation) goes to zero. One considers certain so-called `invertible' systems, which in particular have a unique ground state and a mass gap. It follows from this assumption that all  invariants of the limiting topological theory will necessarily take values in the unit circle. These invariants will be multiplicative under disjoint union of manifolds and will change by conjugation under reversal of orientation. 
Since one is taking limits over finer and finer triangulations, one must assume an independence of the invariants under certain `moves' replacing one triangulation with another. These moves can (and should) be viewed as coming from a triangulation of the product of space-time with the interval by restricting to the two ends of the interval. Since the moves are local, this means that in fact the invariance holds for triangulations of Spin cobordisms between two possibly distinct Spin manifolds, not just triangulations of the product with an interval. 
Thus, such invariants represent elements in the Pontrjagin dual of Spin bordism.\footnote{\label{circle} Recall we identify  the unit circle $S^1$ and $\R / \Z$  via complex exponentiation $e^{2\pi it}$.  Complex conjugation in $S^1$ then identifies with sign change of $t$ in $\R / \Z$.} For a discussion of this in a related case without fermions, see [6].

\subsection{Organization of the Sections of the Paper}

The paper is organized as follows. In the next section, \S1.4, we list some terminology, conventions, and notations that are used throughout the paper.  We conclude the introduction  in \S1.5 with a brief review of some simplicial theory of Postnikov towers.  In \S2 we briefly review Spin bordism in dimensions 2 and 3. The point here is to begin organizing results that eventually lead to using suspensions to compare Spin bordism groups and their Pontrjagin duals in different dimensions. \\ 

In \S3, before turning to technical details and actually constructing the group $G^4(X)$, we  outline the construction of  the dual pairing $G^4(X) \times \widetilde{\Omega}_4^{spin}(X) \to \RR / \ZZ$. It makes sense to do this early on, since certainly two very important points of the paper are describing precisely the cochain triples $(w,p,a)$ that represent elements of $G^4(X)$ and describing how these triples evaluate on Spin bordism representatives $M^4 \to X$.  There are some very delicate details involved in the main result that our pairing does define a group isomorphism $G^4(X) \simeq {\rm Hom}(\widetilde{\Omega}_4^{ spin}(X),\Ar/\Zee)$.  In \S3 we clarify the crucial issues, and indicate where in the paper these various details are dealt with.  In a sense then, \S3 is the most important section of the paper, and the only place where a full proof of the main result is outlined.  The other sections  basically just provide some details of various necessary steps in the proof.  \\

In  \S4 we define quite explicitly the set of cochain triples $(w, p, a)$ that will represent elements of our group $G^4(X)$. Describing these triples is equivalent to giving an explicit cochain description of the 3-stage Postnikov tower that represents the functor $G^4(X)$. The second k-invariant of the classifying Postnikov tower for $G^4(X)$  is  a secondary cohomology operation.  A cocycle formula for this k-invariant involves a $\Z/2$ cochain $x(a)$ that satisfies $$dx(a) = Sq^2Sq^2(a) + Sq^3Sq^1(a) = Sq^2(a^2) + (Sq^1(a))^2$$ for 2-cocycles $a$. Such a cochain $x(a)$ exists because of the Cartan formula for Steenrod squares, or because of the Adem relation $Sq^2Sq^2 + Sq^3Sq^1 = 0$.  In \S5 we establish some important properties of our choice of $x(a)$ that are needed to construct the product of triples and the relations between triples. In \S6 and \S7 we then give details about the product formula for triples and the relations between triples that explicitly lead to the definition of  our group $G^4(X)$.\\

In \S8.1-8.2 we discuss a suspension operation on cochains and apply this to establish an isomorphism between the group $G^3(X)$ of our previous paper [1], which is Pontrjagin dual to $\Omega_3^{spin}(X)$, and the group $G^4(\Sigma X)$ of this paper.  In \S8.3 we use suspensions to fill in a detail that is very important in the dual pairing discussion of \S3. \\ 

In \S9 we deal with a subtle detail that resolves a previous ambiguity from \S3 in the definition of the dual pairing $G^4(X) \times \widetilde{\Omega}_4^{spin}(X) \to \RR / \ZZ$.  This required some computer computations with a complicated  triangulation of $I \times S^2 \times S^2$. The ambiguity itself was analogous to a similar ambiguity that we encountered in our 3-dimensional study [1], which required computations with a complicated triangulation of $\R P^3$.\\

Finally in \S10 we include three appendices.  The first appendix just states some explicit cochain formulas for terms that appear in our formulas for the product of triples in \S6.1 and the relations among triple in \S7.1.  The second appendix provides some information about the group extensions that arise in the Atiyah-Hirzebruch spectral sequence for $4D$ Spin bordism.  These extensions are equivalent to extensions arising from a natural filtration of our group $G^4(X)$, but because our group $G^4(X)$ is constructed from cochains we are able to give formulas that determine these extensions, up to isomorphism. The third appendix studies a canonical subgroup of the Pontrjagin dual of $\widetilde{\Omega}^{spin}_n(X)$, for any $n$.  This subgroup when $n = 4$ is important in the  discussion of the evaluation $G^4(X) \to {\rm Hom}(\widetilde{\Omega}_4^{spin}(X),\Ar/\Zee)$  in \S3.  The subgroup is no harder to study for arbitrary $n$, and also gives rise to an interesting combinatorial characterization of equivalence classes of Spin structures on $n$-manifolds, as certain `quadratic functions' defined on $Z/2$ cocycles of degree $n-1$.  This characterization can also be found in the physics paper [4].\\

\subsection{Terminology, Conventions, and Notations}

Throughout we will work with simplicial sets and with spaces that are {\em ordered} simplicial complexes, that is, the vertices are partially ordered so that the vertices of each simplex are totally ordered.  Maps will be (weakly) order preserving simplicial maps.  Simplicial complex will always mean ordered simplicial complex and simplicial map will always mean ordered simplicial map.\\

We  work with {\em normalized} cochain complexes, consisting of  cochains that vanish on degenerate simplices.  Thus  cochain and cocycle always means normalized cochain and cocycle.  If a simplicial complex $X$ is fixed or understood, we will  write $C^*( F), Z^*(F), H^*(F)$ to indicate  cochains, cocycles, and cohomology of $X$ with coefficients in an abelian group $F$.\\ 

Ordered simplicial structures allow us to compute cup products and higher cup$_i$ products in various cochain complexes.  We use the standard formulas of Steenrod for the cup$_i$ products, and we make use of standard properties of these products, especially the coboundary formula.  For integral cochains $X, Y$ of degrees $|X|, |Y|$ that formula is: 
\begin{eqnarray*}
\lefteqn{d(X \cup_i Y) = }\\
 & &(-1)^i\bigl ( dX \cup_i Y + (-1)^{|X|}X \cup_i dY - X \cup_{i-1} Y - (-1)^{i+|X||Y|} Y \cup_{i-1}X \bigr).
\end{eqnarray*}

The $k$-invariants and product formulas and relations between cochain tuples that we encounter in our Postnikov towers are mostly expressed in terms of cup$_i$ operations.  For computer computations, we turned to classical expressions of cup$_i$ operations as sums of products of evaluations of cochains on faces of simplices.\\

For a $\Z/2$ cochains $c$, we repeatedly encounter the  cochain definition of the Steenrod square operations. That definition is $$Sq^i c =  c \cup _j c + c \cup_{j+1} dc,\ where\ i+j = |c|.$$ If $dc = 0$, this is the standard formula for $Sq^i$ on cocycles.  The cochain formula has the nice property $Sq^i(dc) = dSq^i (c)$.  The shorthand notation $Sq^i c$ doesn't  help with actual computation, but it makes expressions of many of our formulas and proofs of statements about those formulas  more efficient.\\

Many of our cochain formulas involve {\em special lifts} of $\Z/2$ cochains $c$.  The special lift is the $\Z$ cochain $C$ that only takes values 0 and 1 on simplices and reduces to $c \pmod 2$. Special lifts played an important role in our $3D$ work in [1], and are important again in our $4D$ theory.  If $a$ is a 1-cocycle then $dA = 2 A^2$ and if $b$ is a 2-cocycle then $dB = 2B \cup_1 B$.  These formulas can be checked by direct evaluations on simplices.  \\

We frequently encounter functors\footnote{One can make categorical sense of this by working with objects consisting of spaces together with a tuples of cochain, and the obvious notion of map between such objects.} of cochains, $\phi(c_1, c_2, ..., c_k)$, assigning a new cochain  to some cochain or set of cochains, with an obvious naturality property for maps between spaces. We also refer to such cochain functors $\phi$ as {\em natural cochain operations}. Examples include the special lift functor $C = C(c)$, cochain maps induced by coefficient morphisms, and cup$_i$ products.  Sometimes natural cochain operations are only interesting, or only have good properties, when some or all of  the input variables are cocycles.  In our $3D$ work ([1], \S3.2), we discussed and proved some results about cochain functors whose values are cocycles.  Those results are also relevant for  parts of our $4D$ work here as well, especially in \S6 and \S7 in proofs of associativity and commutativity of certain products.\\

Other examples of  natural cochains $\phi$ arise as solutions of `differential equations' $d\phi(c_j) = \Phi(c_j)$, where $\Phi$ is a given natural cochain operation.  For example, $\Phi$ might be a cocycle formula for an Adem relation, which is then a natural coboundary. The existence of  natural solutions $\phi$ to such equations can often be deduced by simply appealing to the existence of a solution in a universal example, such as a product of Eilenberg-MacLane spaces, or a more general Postnikov tower, which we discuss in the next section.  But one can  sometimes just write down a natural cochain $\phi$ that satisfies  a given differential equation, so that proves it exists!  Note the difference of two natural solutions will always be a natural cocycle function of the variables. \\ 

Most of the natural cochain operations that we encounter are expressed as sums of products of evaluations of cochains on faces of simplices. For the purposes of formulas, we write the evaluation of $n$-cochains  $c$ on  ordered $n$-simplices, as $c(012...n)$.  Simplices of general simplicial sets are not  determined by their vertices, but face operators in formulas can still be expressed by simply deleting integers in  sequences $(012...n)$.  As examples, if $p+q = n$, the cup product of a $p$ cochain and a $q$-cochain is given by $$ (a \cup b)(012...n) = a(0...p)b(p...n).$$  If $A$ and $B$ are integral 2-cochains then $$(A \cup_1 B)(0123) = A(013)B(123) - A(023)B(012).$$  It is obvious that such sums of products of face evaluations are functorial in the cochain variables, with respect to simplicial maps between spaces.\\

\subsection{Simplicial Models for  Postnikov Towers} 

In this section, we briefly recall some classical algebraic topology concerning the simplicial theory of Postnikov towers.  By introducing this rather generally here, the special cases that we need below for our study of Spin bordism are easily seen to fit into the general framework.\\

Suppose $A_1$ is an abelian group.  A simplicial set model for an Eilenberg-MacLane space $E_1 =  K(A_1, n_1)  $ has as the $q$-simplices the set of $A_1$-valued normalized $n_1$-cocycles on the standard $q$-simplex $\Delta^q$, with the standard face and degeneracy operations. Note $K(A_1, n_1)$ has a tautological fundamental cocycle  in $Z^{n_1}(K(A_1, n_1); A_1)$.  If $X$ is a simplicial complex, or if $X$ is a simplicial set, then a simplicial map $X \to E_1$ is exactly an $A_1$-valued $n_1$-cocycle, say $a$,   on $X$.  The  null-homotopic maps are those that extend to a simplicial map $CX \to E_1$, where $CX$ is the cone on $X$. This is equivalent to saying $a = dp$ for some $n_1-1$ cochain $p$ on $X$. Thus homotopy classes of simplicial maps $[X, E_1] = H^{n_1}(X; A_1)$. \\

Next, we discuss 2-stage Postnikov towers $$E_2 = K(A_1, n_1) \ltimes_{k(a)} K(A_2, n_2).$$ The notation is meant to indicate a principal fibration over the first Eilenberg-MacLane space with fiber the second Eilenberg-MacLane space.  The term  $k(a)$ is a cocycle representing the cohomology $k$-invariant of the fibration  $$ \bar{k}( a) \in H^{n_2+1}(K(A_1, n_1); A_2).$$ In practice, $k(a)$ is a natural cochain level version of a cohomology operation.  The q-cells of $E_2$ are given by pairs $(p, a)$, where $ a$ is an $A_1$-cocycle on $\Delta^q$ and $p$ is an $A_2$-cochain on $\Delta^q$ with $dp = k(a)$.  Then a simplicial map $X \to E_2$ is given by a pair of cochains $(p, a)$ on $X$, with appropriate coefficients and dimensions, with $da = 0$ and $dp = k(a)$.  Maps $(p_0, a_0)$ and $(p_1, a_1)$ are homotopic if there is a pair of cochains $(\hat{p}, \hat{a})$ on $I \times X$ with $d\hat{a} = 0$ and $d\hat{p} = k(\hat{a})$ that restricts to $(p_i, a_i)$ on $\{i\}\times X$, for $i = 0, 1.$\\

We can continue and define 3-stage Postnikov towers 
$$ E_3 = K(A_1, n_1) \ltimes_{k(a)}  K(A_2, n_2) \ltimes_{k(p, a)} K(A_3, n_3).$$  Thus $E_3$ is a principal  fibration over $E_2$, with an Eilenberg-MacLane space fiber and $k$-invariant determined  by a natural cocycle $k(p, a) \in Z^{n_3+1}(E_2; A_3)$.  Simplicial maps $X \to E_3$ are given by triples of cochains $(w, p, a)$ on $X$, with appropriate coefficients and dimensions, with $da = 0$, $dp = k(a)$, and $dw = k(p, a)$. Homotopic triples are defined just as in the 2-stage case.\\

A tower represents an abelian group valued homotopy functor if it is an H-space, that is, there is a  homotopy associative and commutative simplicial product $E \times E \to E$, with a homotopy inverse. In the 3-stage case, this  means a product of triples $(w, p, a)(v, q, b) = (u, r, c)$, with appropriate properties. In the case of an H-space, the null-homotopic triples determine all relations between triples because of the group structure.  An allowed triple $(w,p,a)$ on $X$ represents the 0 element in $[X, E]$ if there is an allowed triple $(\hat{w}, \hat{p}, \hat{a})$ on the cone $CX$  that restricts to $(w,p,a)$ on the base $ X$. \\

\section {Review of 2- and 3-Dimensional Spin Bordism}

\subsection{The Pontrjagin Dual of $2D$ Spin Bordism} 

The first non-trivial example of Spin bordism is the $2$-dimensional case, since reduced 1-dimensional Spin bordism of $X$ is naturally  isomorphic to $H_1(X, \Z)$.  An element of the group $G^2(X)$ Pontrjagin dual to reduced  $2$-dimensional Spin bordism  is determined by a pair of cocycles $$(p,a) \in  Z^2(X; \R/ \Z) \times  Z^1(X; \Z/2).$$ The only relations are $(p+ dc, a+dx) \equiv (p, a)$, so an equivalence class is just a pair of cohomology classes. The product on the cochain level is given by $$(p,a)(q,b) = (p+q + (1/2)ab, a+b),$$ where $(1/2)$ means the coefficient morphism $\Z/2 \to \R / \Z$. So the classifying space is just $K(\Z / 2, 1) \times K(\R / \Z, 2) $, but the H-space structure is twisted.  Note the cochain product itself is not commutative, but yields an abelian group $G^2(X)$ of equivalence classes because of the cochain formula $ab+ba = d(a \cup_1 b)$.  \\

If $\Sigma$ is an oriented surface,  Spin structures on $\Sigma$ correspond canonically with functions $q\colon H^1(\Sigma ; \Z/2) \to \Z/2$ that  are quadratic refinements of the cup product pairing. We describe the evaluation of a pair of  cocycles $(p,a)$ on a reduced Spin bordism element $f\colon \Sigma^2\to X$. Let $\bar{a}$ denote the cohomology class represented by the cocycle $a$.  Let $Z_a$ be an embedded  $1$-manifold in $\Sigma$ dual to the cohomology class  $f^*\bar{a}$. The Spin structure on $\Sigma$ determines a Spin structure on $Z_a$.  Letting $[Z_a] \in \Omega_1^{spin}(pt) = \Z/2 $ denote  this Spin bordism class, then $q(f^*\bar{a}) = [Z_a] \in \Z /2$.  The evaluation is then given by
$$\langle (p, a), (\Sigma \xrightarrow{f} X) \rangle =    (1/2)q(f^*\bar{a}) +  \int _ \Sigma f^*p \   = \ (1/2)[Z_a]  +  \int _ \Sigma f^*p   \in \R / \Z.$$

It is perhaps amusing that the obvious short exact sequence $$0 \to H^2(X; \R/Z) \to G^2(X) \to H^1(X; \Z/2) \to 0$$ does split for any $X$.  Splittings exist because the isomorphism class of the extension is determined by the function  $$\bar{a} \mapsto (0,\bar{a})^2 = ((1/2)\bar{a}^2, 0) \in H^2(X; \R / \Z) / 2H^2(X; \R / \Z),$$ but $(1/2)\bar{a}^2 = 0$ since $\bar{a}^2$ is the reduction of the integral Bockstein torsion class $\beta \bar{a}$. However, there is no natural group splitting, hence no natural dual pairing of the direct product of cohomology groups and reduced Spin bordism.  Note in the natural pairing above the occurrence of the quadratic function $q$, with $q(\alpha + \beta) = q(\alpha) + q(\beta) + \langle \alpha \beta, [\Sigma] \rangle$, is consistent with the  product formula $(p+q+(1/2)ab,\ a+b)$ in $G^2(X)$.\\

\subsection{The Pontrjagin Dual  of $3D$ Spin Bordism}

Given a simplicial complex $X$, in [1] we constructed a group $ G^3(X)$, functorial in $X$,  from  equivalence classes of  cochain triples on $X$, $$(w,p,a) \in C^3(\RR/ \Z) \times Z^2(\ZZ/2) \times Z^1(\ZZ/2),$$ that satisfy $dw = (1/2)p^2.$   The classifying Postnikov tower is $$K(\Z/2, 1) \times K(\Z/2, 2) \ltimes_{(1/2)p^2}  K(\R/ \Z, 3).$$ We studied a product on triples, given by  $$(w,p,a)(v,q,b) = ( w+v+u,\ p + q + ab,\ a+b)$$ where  $$u = (1/2)p \cup_1q + (1/2)(p+q)\cup_1 ab + (1/2)a(a\cup_1b)b + (1/4)A^2B.$$    Here, $A$ and $B$ denote the special $\Z$ lifts of the $\Z/2$ cocycles $a$ and $b$, which, recall, means the cochain lifts that take only values 0 and 1 on simplices.\footnote{In our paper [1] a term $(1/4)AB^2$ occurs in the product formula for $G^3(X)$ instead of $(1/4)A^2B$.  But in $G^3(X)$ the term $u$ is only well-defined up to coboundaries and since $dA = 2A^2$ we have $(1/4)A^2B - (1/4)AB^2 = d((1/8)AB)$.  So the product formula defined here and the product formula of [1] are actually identical.} The $(1/4)$ in the product formula means the coefficient morphism $\Z \to \R / \Z$ taking 1 to $1/4$.\\

The null-homotopic triples are given by $$(df+(1/2)tdt, dt, dx) \equiv (0, 0, 0).$$ We proved in ([1], \S3) that, for the given product,  triples modulo null-homotopic triples form an abelian group $G^3(X)$.  Note from \S1.4, $tdt = Sq^2t$ for a cochain $t$ of degree 1.\\

We then defined  a natural isomorphism $G^3(X) \simeq \rm{Hom}(\Omega_3^{spin}(X), \RR/\ZZ)$.  From the product formula, representatives of elements in $G^3$ factor as $$(w,p,a) = (w,p,0)(0,0,a).$$  Thus the evaluation of $(w,p,a)$ on a Spin bordism element must be the sum of the evaluations of $(w,p,0)$ and $(0,0,a)$.\\

The evaluation of elements $(w,p,0)$ on a Spin bordism class $f\colon M^3 \to X$ was described in two equivalent ways.  Details can be found in ([1], \S6), and are quite similar to the $4D$ analogues in \S3.4 below. The evaluation of an element $(0,0,a)$ was defined to be the Arf invariant\footnote{\label{Arf}The natural values of Arf invariants are $8^{th}$ roots of unity in the unit circle. Using the complex  exponential identification $\R/ \Z \simeq S^1$, we identify Arf invariants with elements of the subgroup  of $\R /\Z$, generated by $1/8$.}  in $(1/8)\Z/ \Z \subset \R / \Z$ of a ${\rm Pin}^-$ surface $\Sigma \subset M^3$ dual to the cohomology class $f^*(\bar{a}) \in H^1(M^3; \Z/2)$.\\

The proof that our evaluation was a group homomorphism from $G^3(X)$ to the dual of Spin bordism required some care with Arf invariants.  The manifold $\R P^3$, given a Spin structure as the boundary of $T_DS^2$, the tangent disk bundle of $S^2$ with its natural orientation, contains a Pin$^{-}$ surface $\R P^2$ dual to the fundamental class  $\bar{a} \in H^1(\R P^3; \Z/2)$.  The standard convention is that $\rm{Arf}(\R P^2) = +1/8$, which is equivalent to the underlying quadratic function  $q\colon H^1(\RR P^2; \ZZ/2) \to \ZZ/4$ satisfying $q(\bar{a}) = +1$.  In ([1], \S8) we showed that the product on $G^3(X)$ defined above was consistent with this $\rm{Arf}(\RR P^2) = +1/8$ choice in the evaluation formula.  Essentially this amounted to evaluating $(0,0,a)^2 = ((-1/4)A^3, a^2, 0) \in G^3(\R P^3)$ on the identity map $\R P^3 \to \R P^3$, and getting the result $+1/4$, using a rather complicated triangulation of $\R P^3$.\\ 

At this point it is appropriate to mention two details about the $3D$ theory that end up complicating the $4D$ theory.  The first is that there is a second product possible for $3D$ triples, which replaces the term $(1/4)A^2B$ in the product formula above by $(-1/4)A^2B$.  The result is an isomorphic group, but if we used that product in $G^3(X)$ we would have changed the evaluation of elements $(0,0,a)$ on $f \colon M^3 \to X$  to be $(-1)\rm{Arf}(\Sigma)$, where $\Sigma \subset M$ is a Pin$^{-}$ surface dual to $f^*(\bar{a})$.   But, more to the point, with our choice of product in $G^3(X)$, a certain choice of $k$-invariant is forced in the $4D$ theory in order that our suspension map $s\colon G^3(X) \to G^4(\Sigma X)$ of \S8.2 is an isomorphism.  We discuss this issue further in later sections. \\

The second detail is that with our product fixed on $G^3(X)$ there is an automorphism, different from the inverse automorphism, given by $(w,p,a) \mapsto (w+(1/2)a^3, p, a)$.\footnote{So the automorphism group of the functor $G^3(X)$, hence also of $\Omega_3^{spin}(X)$, is $\Z/2 \oplus \Z/2 = Aut(\Z/8)$.  It is amusing that this is  easier to see with $G^3(X)$ than with $\Omega_3^{spin}(X)$.}  This can be seen to imply that in the evaluation formula for our choice of $3D$ product, we could have evaluated $(0,0,a)$ as either $\rm{Arf}(\Sigma) $ or $5\rm{Arf}(\Sigma)$ in $(1/8)\Z /\Z$.   The significance of this second detail for the $4D$ theory is discussed in \S3.1 and \S3.3 below. \\  

\section{Preview of the Dual of $4D$ Spin Bordism}

\subsection{Outline of the Isomorphism $G^4(X) \simeq {\rm Hom}(\widetilde{\Omega}_4^{ spin}(X),\Ar/\Zee)$}

In this section we will outline  the details involved in constructing the group $G^4(X)$ and identifying that group with the Pontrjagin dual of reduced $4D$ Spin bordism. In fact, this is the only place in the paper where the full proof is organized.  All other sections deal with particular details. So this outline is crucial for comprehending the overall proof.\\

 Instead of viewing the  groups $G^n(X)$ Pontrjagin dual to Spin bordism in different dimensions as independent constructions, they should be viewed as groups in a generalized cohomology theory, represented by an $\Omega$-spectrum of Postnikov towers, with suspension isomorphisms $$s \colon G^n( X) \to G^{n+1}(\Sigma X).$$  The separate evaluation isomorphisms are chosen so that they commute with the Pontrjagin duals of  geometrically defined bordism suspension isomorphisms $$\widetilde{\Omega}_{n+1}^{spin}(\Sigma X) \to \widetilde{\Omega}_n^{spin}(X).$$ For example, once we clarify some properties of  the cochain suspension maps $s$ in \S8.1, it is fairly easy to prove that the suspension $s \colon G^2(X) \to G^3(\Sigma X)$ defined by $s(p, a) = (sp, sa, 0)$ is an isomorphism, and commutes with the evaluations on Spin bordism defined in \S2.\\

We will now give a very brief summary of our preview of the $4D$ theory, which we will expand below.  First we write down a 3-stage Postnikov tower $E$ that represents the Pontrjagin dual of $4D$ Spin bordism. Thus $G^4(X) = [X, E]$, and elements of $G^4(X)$ are equivalence classes of certain triples of cochains. There are two cohomology classes possible for the second $k$-invariant  of $E$, differing by sign.  Each choice leads to  two products of $4D$ triples $(w, p, a)$.  But for only one of these $k$-invariant choices are the $4D$ products suspension compatible with the $3D$ product of \S2.2.  So that determines our $k$-invariant choice for $E$.  Then we turn to the evaluation isomorphism $$G^4(X) \simeq {\rm Hom}(\widetilde{\Omega}_4^{\rm spin}(X),\Ar/\Zee).$$  It turns out that for either choice of product on $G^4(X)$ there is a unique  evaluation isomorphism, compatible with the simple evaluation formula $$\langle (w, 0, 0), (M \xrightarrow{f} X) \rangle = \int_M\ f^*w.$$  In order to be suspension compatible with a $3D$ evaluation, a few other necessary $4D$ evaluation formulas must hold.  There was a choice of two possible $3D$ evaluations, mentioned at the end of \S2.2.  However, it turns out that  the $3D$ evaluation we actually defined in [1], and in \S2.2 above, can  be extended to a $4D$ evaluation that is additive for only one of the two possible $4D$ products. So that is the product we choose on $G^4(X)$ in \S6. Identifying which $3D$ evaluation extends to which $4D$ product and evaluation was one of the most difficult parts of our work.  So that's it.  Now we elaborate.\\

The details of the construction of the group $G^4(X)$ and the $4D$ evaluation isomorphism are rather complicated.  First, we need to specify a cochain level version of the 3-stage Postnikov tower that represents the Pontrjain dual of $4D$ Spin bordism.  That is carried out in \S4.   In the language of \S1.5, the $4D$ tower will have the form
$$ E = K(\Z/2, 2) \ltimes_{k(a)}  K(\Z/2, 3) \ltimes_{k(p, a)} K(\R / \Z, 4).$$
The homotopy groups are just the duals of reduced Spin bordism groups of spheres. Simplicial maps $X \to E$ are named by triples of cochains $(w,p,a)$ with $dp = k(a)$ and $dw = k(p,a)$.\\

The first $k$-invariant is cohomologically $Sq^2\bar{a} = \bar a^2$.  The cohomology class of the second $k$-invariant lies in a group isomorphic to $\Z/4$, and has order 4. One could choose either generator, since the two choices give homotopy equivalent towers.  We will discuss in \S4.3 why these towers do indeed represent the dual of reduced $4D$ Spin bordism.  With either choice, the loop space is the classifying tower for the dual of $3D$ Spin bordism, discussed in \S2.2.   However, a specific  choice  of the second $k$-invariant influences cochain formulas  for the product of allowable triples in $ G^4(X)$.  In \S8.2 we show that only one of the cohomology $k$-invariant choices will yield a group $G^4(X)$ for which the suspension $s\colon G^3(X) \to G^4(\Sigma X)$ defined by $s(\omega, \rho, \alpha) = (s\omega, s\rho, s\alpha) = (w,p,a)$ is an isomorphism.  So that is the $k$-invariant we choose. If we had chosen the alternate product in $G^3(X)$ in [1], we would have chosen the alternate $k$-invariant for the $4D$ theory here.\\   

Next, once the $k$-invariant is fixed, in \S6 we study product formulas for triples $(w,p,a)(v,q,b) \in G^4(X)$.  There are  two possible product formulas in the $4D$ case, discussed in \S6.4, leading to isomorphic groups and both desuspending to the product in $G^3(X)$. \\

Finally we look at possible evaluation formulas. In the  $4D$ case, our evaluation homomorphism $ G^4(X) \to \rm{Hom}(\widetilde{\Omega}^{spin}_4(X), \RR/\ZZ)$ will also be divided into two  cases, similar to the $3D$ case.  However, the reduction to the two cases is harder since elements  $(w,p,a) \in G^4(X)$ do not easily factor. Instead, given $f\colon M^4 \to X$, then in \S8.3 we show elements  $(w, p, a) \in G^4(X)$ pulled back to $G^4(M)$ do factor, $(f^*w, f^*p, f^*a) = (w',p',0)(0,0,a') \in G^4(M)$, after possibly subdividing the triangulation of $M$. By naturality, we must have equality of evaluations $$\langle (w, p, a),\ (M \xrightarrow{f} X) \rangle = \langle (f^*w, f^*p, f^*a), (M \xrightarrow{Id} M) \rangle$$ $$= \langle (w', p', 0), (M \xrightarrow{Id} M) \rangle + \langle (0, 0, a'), (M \xrightarrow{Id} M) \rangle.$$  There are two equivalent ways to define the evaluation of elements $(w', p', 0)$ on $Id \colon M \to M$, analogous to the $3D$ case, which we discuss in \S3.4. Elements $(0, 0, a')$ again evaluate in terms of Arf invariants, as we discuss in \S3.3.  But there are two possibilities, which turn out to correspond to the two possible products on $G^4(X)$. Determining exactly which Arf invariant choice corresponded to  which $4D$ product formula was difficult and is carried out in \S\S9.2-9.3.  The evaluation of $(w,p,a)$ on $f\colon M \to X$ must then be the sum of these separate evaluations.\\

A priori, it would seem that we need to prove our evaluation is independent of choices made, such as the factorization in $G^4(M)$, and that it is additive in the $G^4(X)$ variable. That would be hard. But, at this point we know the group $G^4(X)$ is isomorphic to the Pontrjagin dual of reduced $4D$ Spin bordism, because we have used a correct Postnikov tower $E$ that represents the Pontrjagin dual and  constructed $G^4(X)$ so that  $G^4(X) = [X, E]$. Moreover, an isomorphism between $G^4$ and the Pontrjagin dual will exist preserving natural filtrations. The filtrations are determined on the bordism side by the skeleton filtration of $X$  and on the $G^4(X)$ side by sub towers and quotient towers of the Postnikov tower $E$.\\

It can be seen that the functor $G^4(X)$ admits no natural automorphisms preserving the filtration, other than the inverse map of abelian groups. This is explained in \S9.1. Therefore, for each possible product on $G^4(X)$, there exists a unique natural evaluation group isomorphism $ G^4(X) \to \rm{Hom}(\widetilde{\Omega}^{spin}_4(X), \RR/\ZZ)$ preserving filtrations, up to a sign choice.  The sign choice can be specified by requiring $$\langle (w, 0, 0), (M \xrightarrow{f} X) \rangle = \int_M\ f^*w,$$ which is part of our definition of evaluation of elements $(w, p, 0)$ in \S3.4. Then, the evaluation  $\langle (0, 0, a'), (M \xrightarrow{Id} M) \rangle $ turns out to be forced to be 1 or 5 times an easily defined Arf invariant.  We discuss this in \S3.3 below. The question is, which  multiple goes with which $4D$ product?\\    

It turns out that we can identify which choice of $4D$ product corresponds to the odd multiple 1 of the Arf invariant that we discuss in \S3.3.  This is the choice that gives suspension compatibility with our previously chosen $3D$ evaluation in [1]. So that is the $4D$ product that we choose in \S6.  The proof of the compatibility of the evaluation we define in \S3.3 and our choice of $G^4(X)$ product in \S6 requires study of the special case $X = S^2 \times S^2$ and a complicated ordered triangulation of $I \times S^2 \times S^2$, and this turned out to require extensive computer computations.  The details of the needed $G^4(S^2 \times S^2)$ computations are given in \S9.  \\

\subsection{Preview of a Factorization of Elements of $G^4(M)$}

Given a Spin bordism representative $f\colon M \to X$, we need to define an invariant $$\langle (w,p, a),\ (M \xrightarrow{f} X) \rangle \in \RR / \ZZ.$$  We outlined in the previous section, and will fill in details here and in the next two sections, how that definition goes. Much of this discussion is conceptually easier to assimilate than the technical details of the $G^4(X)$ construction.  Also, the evaluation discussion is likely of greater interest to physicists than the details of the group construction. \\

By the naturally requirement, we must have $$\langle (w, p, a),\ (M \xrightarrow{f} X) \rangle = \langle (f^*w, f^*p, f^*a), (M \xrightarrow{Id} M) \rangle,$$ so we can just assume $X = M^4$ is a Spin 4-manifold and $f$ is the identity.  Now it will turn out that after possibly subdividing the triangulation of $M$, any representative $(w, p, a)$ of an element of $ G^4(M)$ can be multiplied by a relation triple (representing the zero element of the group) and then factored at the level of the product of cochain triples, so that we have in $G^4(M)$ $$(w,p,a) \equiv (w, p, a)(relation) = (w', p', a') = (w', p', 0)(0,0, a').$$  In fact, referring to this factorization of representatives on the right, there will exist a triangulation of $M$ and a simplicial map $g\colon M \to \Sigma \RR P^\infty$ so that $ a' = g^*(s\alpha)\in Z^2(M; \ZZ/2)$ is the pull back by $g$ of a  cocycle $s\alpha \in Z^2(\Sigma \RR P^\infty; \ZZ/2)$ with $(s\alpha)^2 = 0$.  Specifically, $\alpha \in Z^1( \R P^\infty; \Z /2)$ is a cocycle that generates the cohomology, and $s$ is a cochain suspension map defined and studied in \S8.1 below.  Moreover, $(0, 0, s\alpha)$ represents an element of $G^4(\Sigma \R P^\infty)$, and then $g^*(0, 0, s\alpha) = (0, 0, a') \in G^4(M)$.\\

This factorization  $(w,p,a) = (w', p', 0)(0, 0, a') \in G^4(M)$ will be discussed further in \S8.3 below.  We emphasize that the explicit formula for $(w', p', 0)$ in terms of $(w, p, a)$ and other cochains is very complicated.  An explicit formula is essentially given in \S8.3. Be that as it may, we see that we must have $$\langle (w,p, a),\ (M \xrightarrow{Id} M) \rangle  =$$  $$ \langle (w',p', 0),\ (M \xrightarrow{Id} M) \rangle + \langle (0,0,s\alpha), (M \xrightarrow{g} \Sigma \RR P^\infty) \rangle.$$  We will separately define the two summands on the right side of this last equation  in the next two subsections.  \\

\subsection{The Evaluation $\langle (0, 0, s\alpha), (M \xrightarrow{g} \Sigma \R P^\infty) \rangle $}  

There is a geometric suspension isomorphism $\widetilde{\Omega}_4^{spin}(\Sigma X) \simeq \widetilde{\Omega}_3^{spin}(X)$, defined by making a map $g\colon M^4 \to \Sigma X$ transverse to $X$, resulting in a collared Spin 3-manifold $(-1, 1) \times N^3  \subset M^4$ and a map $g\colon N^3 \to X$. Here we interpret the suspension $\Sigma{X}$ as the obvious quotient space of $[-1, 1] \times X$. So the suspension is the union of two cones, $C^+X$ and $C^-X$.\\

The geometric suspension isomorphism requires an orientation choice for $N$.  The natural choice would be the orientation so that the product orientation on $(-1,1) \times N$ agrees with the orientation of $M$.  But  we will choose in this dimension the opposite orientation of $N$.  In fact, the Spin structure that we choose on $N$ can be described as the induced Spin structure on $\partial M_+$, where $M_+ =  g^{-1}(C^+X)$ is oriented with the outward (downward pointing) normal first.  The reason for this choice is that  we will want a commutative diagram of isomorphisms:

\begin{align}\label{diag1}
\begin{split}
G^3(X) \hspace{.3in}&\xrightarrow{s}\ \ \    G^4(\Sigma X)   \\ 
\downarrow \hspace{.5in}& \hspace{.5in}  \downarrow \\
{\rm Hom}(\widetilde{\Omega}_3^{\rm spin}(X),\Ar/\Zee) & \rightarrow  {\rm Hom}(\widetilde{\Omega}_4^{spin}(\Sigma X, \R / \Z) \\
\end{split}
\end{align}

\noindent where $G^3(X)$ is the group we constructed in [1] Pontrjagin dual to $\widetilde{\Omega}_3^{spin}(X)$.\\

There are two ways to see a suspension isomorphism $G^3(X) \to G^4( \Sigma X)$  that will ultimately yield such a  commutative diagram.  A specific group homomorphism $s(v, q, b) = (sv, sq, sb)$ is described in \S8.2 below.  A relatively simple filtration argument proves that this map is an isomorphism.\footnote{A second method is to exploit the fact that our construction of $G^4(Y)$ comes along with a group identification with homotopy classes of maps $[Y,  E]$, where $E$ is the Postnikov tower described in \S4.3.  Therefore, if $Y = \Sigma X$ is a suspension we will have  $[\Sigma X, E] = [X, \Omega E]$.  But the loop space Postnikov tower $\Omega E$ is easily seen to represent the functor $G^3(X)$ of our three dimensional study.} Our choice of $s$ and our choice of the $3D$ evaluation homomorphism in \S2.2 and the $4D$ evaluation here and in \S3.4 dictates our choice of Spin structure on $N$ above, so that Diagram (3.1) will commute on all of $ G^3(X)$.\\ 

We can take $X = \R P^\infty$ and $s(0, 0, \alpha) = (0,0,s\alpha)$. It is well-known, and  is an easy consequence of our general work in [1], that $\widetilde{\Omega}_3^{spin}(\R P^\infty)  \simeq \Z / 8$, with generator $(0, 0, \alpha)$.  We then have $$G^4(\Sigma \R P^\infty) \simeq G^3(\R P^\infty) = \rm{Hom}(\widetilde{\Omega}_3^{spin}( \R P^\infty), \R / \Z) \simeq (1/8)\Z/ \Z.$$  Given a Spin bordism element $g\colon N^3 \to \R P^\infty$, from ([1], \S7.1) the assigned invariant  in $(1/8)\Z / \Z$  is  the Arf invariant of a quadratic function  $q\colon H^1(\Sigma; \Z/2) \to \Z /4$, where $\Sigma \subset N$ is a Pin$^{-}$ surface dual to  $g^*(\bar{\alpha}) \in H^1(N; \Z /2)$. We then define $$\langle (0, 0, s\alpha), (M^4 \xrightarrow{g} \Sigma \R P^\infty) \rangle = \langle (0,0,\alpha), (N^3 \xrightarrow{g} \R P^\infty) \rangle = \rm{Arf}(\Sigma, q) \in \R / \Z,$$ where $M^4$ and $N^3$ are related by the general geometric suspension discussion above.\\

Now, there is a  difficulty here.  The group $G^3(X)$ admits an automorphism $(w,p,a) \mapsto (w+(1/2)a^3, p, a)$.  This means that another possible $3D$ evaluation would be to replace $\rm{Arf}(\Sigma, q)$ in the above formula by $\rm{Arf}(\Sigma, q) +(1/2)\int_N g^*\alpha^3 = 5\rm{Arf}(\Sigma, q)$. But the group structure on $G^4(X)$  admits no natural automorphisms other than the inverse map. This means that the two possible evaluation choices for $(0,0,g^*s\alpha) \in G^4(M)$ correspond to the two possible product choices on $G^4(X)$.\\

It turns out that deciding which $G^4(X)$ product choice corresponds to the $\rm{Arf}(\Sigma, q)$ evaluation choice just declared can be decided with the special case $X = S^2 \times S^2$, and a study of the product $(0,0,a)(0,0,b)$, where $a, b \in Z^2(S^2 \times S^2; \Z/2)$ are pullbacks of a standard cocycle in $Z^2(S^2; \Z/2)$ under the two obvious projections.  Fortunately, it was not necessary to find an explicit factorization of the form discussed above, $(0,0,a)(0,0,b) \equiv (w', p', 0)(0, 0, a+b+dx)$, and then carry out the separate evaluations of the terms on the right side.  This would have required computer computations with the second 
barycentric subdivision of $S^2 \times S^2$, which is too big.  Instead, we were able to work with a somewhat exotic vertex ordering on $I \times (S^2 \times S^2)$, using only  the first barycentric subdivision of $S^2 \times S^2$ .  This change of viewpoint about relations is explained in \S7.1.  The final result, proved in \S9, is that the $\rm{Arf}(\Sigma, q)$ choice in the evaluation formula above is consistent with the $4D$ product formula we choose in \S6.\\

There is a final important point to be made about the Arf invariant evaluation of this section, which is needed in the final filtration arguments in \S7.5 that are used to prove our full evaluation $G^4(X) \to {\rm Hom}(\widetilde{\Omega}_4^{ spin}(X),\Ar/\Zee)$ is an isomorphism.  Suppose $g\colon M^4 \to S^2 \subset \Sigma \R P^\infty$. Then as above we find a framed oriented 3-manifold $N^3 \subset M^4$ and then a framed oriented surface $\Sigma \subset N^3$.  So $\Sigma$ is framed in $M^4$, hence inherits a Spin structure.  In this case, ${\rm  Arf}(\Sigma, q)$ is an element in $(1/2)\Z / Z \simeq \Z/2$, which names the Spin bordism class in $\widetilde{\Omega}^{spin}_2(pt) \simeq \Z/2$ of the framed surface $\Sigma \subset M$.\\

\subsection{The Evaluation $ \langle (w, p, 0), (M \xrightarrow{Id} M)\rangle $ }

For any space $X$ the equivalence classes of triples $$\{(w, p, 0)\ |\ dp = 0,\ dw = (1/2)p \cup_1 p = (1/2)Sq^2p\}$$ will form a subgroup $\tilde{G}^4(X) \subset G^4(X)$.  The restricted multiplication turns out to be quite simple, and is given by $$(w,p,0)(v,q,0) = (w+v+(1/2)p \cup_2q,\ p+q,\ 0).$$  The relations are also easy to describe, and are given by $$(w, p, 0) \equiv  (w, p, 0)(df + (1/2)Sq^2c,\ dc,\ 0)$$ $$= (w + df +  (1/2)Sq^2c + (1/2)p \cup_2 dc,\ p+dc,\ 0).$$  Here we use the standard definition of the operation $Sq^2$ on 2-cochains from \S1.4, which is $Sq^2c = c^2 + c \cup_1 dc$.\\

Given a Spin bordism representative $f\colon M^4 \to X$, we will define in two ways the evaluation $\langle (w, p, 0), (M \xrightarrow{f} X) \rangle \in \R /\Z$.  By naturality, it suffices to define the evaluation when $X = M$ and $f = Id$. The two definitions are direct analogues of two definitions in  our three dimensional work  ([1], \S6.1) for evaluation of triples of this simplified form. In fact, these two evaluations of certain pairs $(w,p)$ extend to Spin bordism in all dimensions.  The general case is explained in Appendix \S10.3.\footnote{In fact, in \S10.3 an equivalence class of Spin structure on $M^n$ is identified with a certain `quadratic function' $Q\colon Z^{n-1}(M; \Z/2) \to \Z/2$, and the evaluation of $(w,p)$ on $Id\colon M \to M$ can be very neatly described as $(1/2)Q(p) + \int_{[M]} w$.}\\

The first definition is due to Kapustin.    The 4-dimensional reduced Spin bordism of both $K(\Z/2, 3)$ and $K(\Z/2, 3) \times K(\Z/2, 3)$ vanishes.  This follows trivially from the Atiyah-Hirzebruch spectral sequence.  Therefore, given $p, q \in Z^3(M; \Z/2)$, we can choose a Spin 5-manifold $W$ with $\partial W = M$ and cocycles $\tilde{p}, \tilde{q} \in Z^3(W; \Z/2) $ extending $p, q$.  We define $$\langle (w, p, 0), (M \xrightarrow{Id} M) \rangle =  \int_{[M]}\ w  + (1/2) \int_{[W]}\ \tilde p \cup_1 \tilde p\  \in\  \R / \Z.$$  Here, $\tilde p\ \cup_1\ \tilde p$ is the standard cochain representative of the cohomology operation $Sq^2 \tilde p$ on cocycles.  By applying this construction involving 5-manifolds three times,  to $p, q$ and $p+q$, and using the fact that $Sq^2 $ vanishes on $H^3(\widehat{W}; \Z/2)$ for closed Spin 5-manifolds  and on $H^3(W, \partial W; \Z/2)$ for Spin manifolds with boundary, it is fairly easy to prove that this evaluation does not depend on the choices made and gives a well-defined homomorphism $$\tilde{G}^4(X) \to \rm{Hom}(\widetilde{\Omega}_4^{spin}(X), \R / \Z).$$
We give more details in \S10.3. \\

 For our second evaluation method for elements $(w,p,0) \in G^4(M)$, one first finds an ordered simplicial map $u \colon M \to S^3$ so that  $p + dc = u^*z$, where $z$ is a fundamental cocycle on $S^3$ non-zero on exactly one 3-simplex, and $c$ is some 2-cochain on $M$.  This may require some subdivision of the triangulation of $M$. Then the inverse image under $u$ of a point on the sphere is a framed 1-manifold $Z \subset M$, which inherits a Spin structure from the Spin structure on $M$.  Denote this Spin bordism class by $[Z] \in \Omega_1^{spin}(pt) \simeq  \Z/2$. Recall in $\tilde{G}^4(M)$ we have $$(w,p,0) \equiv (w + (1/2)Sq^2c+ (1/2)p \cup_2 dc ,\  p+dc,\ 0).$$  We then define $$\langle (w,p,0), (M \xrightarrow{Id} M) \rangle = $$  $$\int_{[M]} w + (1/2)\Bigl( [Z] + \int_{[M]} Sq^2c + \int_{[M]}p \cup_2 dc \Bigr)\in \R / \Z.$$ 
In \S10.3 we give a proof that the two evaluations here for elements $(w, p, 0)$ agree, and define an isomorphism of groups $$\tilde{G}^4(X) \simeq \rm{Hom}(\widetilde{\Omega}_4^{spin}(X) / F_2, \R / \Z),$$ where $F_2 = Image(\widetilde{\Omega}_4^{spin}(X^{(2)}) \to \widetilde{\Omega}_4^{spin}(X))$, with $X^{(2)}$ denoting the 2-skeleton of $X$. Moreover, a Stokes Theorem argument, along with a simple observation about a framed 1-submanifold $Z$ in an $N^3 \subset M^4$, related as in the geometric suspension discussion, proves this isomorphism is consistent with the commutativity of Diagram (3.1) above. The \S10.3 proof that the two evaluations agree seems easier than the methods used in ([1], \S\S6.2-6.4) of our three dimensional work proving that  two analogous evaluations of elements of the form $(w,p,0)$ agree. \\

We are now ready to turn to the actual construction of the group $G^4(X)$, which takes up the next four sections.

\section{Triples of Cochains and the Basic $4D$ Equation}

\subsection{Triples of  Cochains and Compact Topologies}

Give a simplicial complex $X$, we begin with the following set of cochain triples.
$$\widehat{\C} =\widehat{\C}(X) =  \{(w, p, a) \in  C^4(X; \RR/\ZZ) \times C^3(X; \ZZ/2) \times Z^2(X; \ZZ/2)\ \vert\ dp = a^2\}.$$

Cochains of any space with coefficients in a compact abelian group form a compact abelian group, specifically just a  direct product of  copies of the coefficient group, indexed by simplices.  All our operations on sets of tuples of cochains below, such as constructions of certain subsets and subquotient sets of $\widehat{\C}$, are compatible with these compact topologies.  Thus it is pretty routine to see why the group $G^4(X)$ we eventually construct from these triples of cochains has a compact topology.\\

 Certainly any thorough discussion of Pontrjagin duality for groups should include a topological discussion.  But then it is also pretty routine to see why the map we construct $G^4(X) \to \rm{Hom}(\widetilde{\Omega}_4^{spin}(X), \RR /\ZZ)$ is continuous in the compact topologies.  So we will suppress this part of the developement.  In any case, for the sort of spaces $X$ that are of interest, such as finite complexes or classifying spaces $B\Gamma$, the topologies on the groups are pretty simple, consisting of finite groups and maybe some $\RR / \ZZ$ summands, corresponding to $\ZZ$ summands of $H_4(X; \ZZ)$.\\

\subsection{The Basic Equation}

The basic equation, which will restrict the triples in $\widehat{\C}$ that we want to work with,  is $\D(w,p,a) = 0$, where $dp = a^2$ and  $\D\colon \widehat{\C} \to C^5(\RR / \ZZ)$ is defined by
 $$\D(w,p,a) = dw - \big[ (1/2)Sq^2p + (1/4)A (A \cup_1A) + (1/2)x(a) \big].  $$
Here $Sq^2p  = p \cup_1 p + p \cup_1 dp$ is the formula for the operation $Sq^2$ on 3-cochains.  Also, $A$ is the special integral lift of $a$, taking only values 0 and 1 on simplices, and  $(1/n)$ means  the coefficient morphism $\ZZ/n \to \RR /\ZZ$ or $\ZZ \to \R /\Z$.  \\

In the formula for $\D$, the term $x(a) \in C^5(\ZZ/2)$ is a functor of cocycles $a \in Z^2( \Z/2)$,  with  $$dx(a) = a^2 \cup_2 a^2 + (a \cup_1 a)^2 = Sq^2(a^2) + (Sq^1(a))^2.$$
Functors of cocycles were discussed in \S1.4. There is such a natural  $x(a)$ because of the Cartan formula for $Sq^2(\bar a^2)$. The Cartan formula can also be interpreted as the Adem relation $0 = Sq^2Sq^2 + Sq^3Sq^1$ applied to  degree 2 cohomology classes $\bar a \in H^2(\Z/2)$.  Then one can find a universal $x(a)$ by finding $x(a)$ for the tautologous 2-cocycle $a \in Z^2(K(\Z/2, 2) ; \Z/2)$.\\

The standard universal arguments show that any two such functorial $x(a)$ differ by a functorial exact class (coboundary)
plus a linear combination of $aSq^1a$ and $Sq^2Sq^1a$.  To make the basic equation  unambiguous, we need to pin down the  cochain $(1/2)x(a)$.  We will do this in \S4.5 below.\\

An important property of special lifts for 2-cocycles $a$ is $dA = 2A \cup_1 A $. If $dx(a)$ satisfies the above  equation, then the formulas $dp= a^2$ and  $dSq^2p = Sq^2 dp$ and $dA = 2A \cup_1 A$ imply that the following expression is a natural  $\R / \Z$-cocycle $$k(p, a) = (1/2)Sq^2p + (1/4)A (A \cup_1A) + (1/2)x(a).$$ This cocycle  represents the secondary cohomology operation $$\langle \frac{1}{2}Sq^2, Sq^2 \rangle \bar{a},$$
and this operation is the $d_3$ differential of $\bar{a}\in H^2(\ZZ/2)$ in the Atiyah-Hirzebruch spectral sequence for the Pontrjagin dual of reduced four dimensional spin bordism.   The $d_3$ differential is defined on $\bar{a}$ with $ Sq^2(\bar{a})  = \bar{a}^2 = 0$, which is the $d_2$ differential of $\bar{a}$ in the spectral sequence.  Thus the basic equation states that $d_3(\bar{a}) = 0$, since the basic equation writes the secondary operation as a coboundary.\\

Our group $G^4(X)$ will be constructed in \S6 and \S7 below by defining a product of two triples $(w,p,a)$ and $(v,q,b)$  that satisfy the basic equation and also defining an equivalence relation between triples, so that the equivalence classes with the induced product form an abelian group.\\

\subsection{The Postnikov Tower}

The homotopy functor that we are studying, the Pontrjagin dual of reduced $4D$ spin bordism, is represented by the following 3-stage Postnikov tower: 
$$E = K(\Z/2, 2) \ltimes_{a^2} K(\Z/2, 3) \ltimes_{k(p,a)} K(\R/\Z, 4).$$
Pairs $(p,a)$ on an ordered simplicial complex $X$ with $dp = a^2$ correspond to   simplicial maps $$X \to K(\Z/2, 2) \ltimes_{Sq^2a} K(\Z/2, 3).$$
The triples $(w,p,a) \in kernel(\D)$, that is, with $dw = k(p,a)$,  correspond to simplicial maps $X \to E$.  Below we will define a product of triples, that is, an H-space structure on $E$, and we will describe the null-homotopic triples, or equivalently  the relations between triples.\\

The second $k$-invariant of the Postnikov tower $E$, $$k(p,a) = (1/2)Sq^2p + (1/4)A (A \cup_1A) + (1/2)x(a),$$ represents a class in the cohomology group $$H^5(K(\Z/2, 2) \ltimes_{a^2} K(\Z/2, 3) ; \R/\Z) \simeq \Z / 4.$$  In fact, $k(p,a)$ represents an element of order 4.\\

We address the question, why does $E$  represent the Pontrjagin dual of reduced $4D$ spin bordism?  First, the homotopy groups are correct.  Second, the loop space of the correct tower $E$, with its loop space product, must represent the Pontrjagin dual of $\Omega_3^{spin}(X)$, and we understand that representing $H$-space from our work in [1], outlined in \S2.2 above.  In particular, we must have $\Z /8 \simeq [\R P^\infty, \Omega E] \simeq [\Sigma \R P^\infty, E]$ under an $H$-space product on $E$ consistent with the $k$-invariants. This is enough to imply that the correct first $k$-invariant of $E$ cannot be 0 and also that the second $k$-invariant of $E$ cannot loop down to 0.  So the second $k$-invariant of $E$ cannot be 0 or the element of order 2, since the element of order 2 is $(1/2)aSq^1a$ which does loop down to 0. Finally, the two second $k$-invariants of order 4 do give homotopy equivalent towers. So we could choose either one of these.\\

Because we work with cochains and cocycles, not cohomology classes, our formulas for the $k$-invariant and for the product of triples and for the relations between triples are highly non-unique, in the sense that many different formulas lead ultimately to isomorphic groups of equivalence classes of triples.  This non-uniqueness is discussed further  in \S4.4 and \S7.6 below.\\  

\subsection{Variations of the Basic Equation}

There are several ways to vary the basic equation formula  $dw = k(p,a)$, so that triples satisfying one equation match up bijectively with triples satisfying the other equation.  For example, one can add any natural coboundary $dc(p, a)$ to $k(p,a)$.\footnote{The term natural coboundary should include or imply $c(0,0) = 0$.} Then  $(w,p,a)$ would be matched with $(w+c(p,a), p, a)$.  Coboundaries can be found that reverse the order of the cup and cup$_i$ product terms in the formula.  Such a change in the basic equation would dramatically alter the product formula for allowable triples, and the relations between triples that we will study below. Next, one can vary the choice of $x(a)$ by adding a linear combination of $aSq^1a$ and $Sq^2Sq^1a$.  Third, there is a free choice of sign $$(\pm 1/4)A(A \cup_1A).$$
 Adding a coboundary  to $k(p, a)$ does not change the cohomological $k$-invariant of the Postnikov tower.  The sign change replaces the $k$-invariant by its negative, which does not change the homotopy type of the Postnikov tower.\\

Regarding changes of $x(a)$,  it turns out that $(1/2)(Sq^2Sq^1a + aSq^1a)$ vanishes in $\RR/\Zee$-cohomology.\footnote{\label{x(a)}Proof: Note $d(A^2) = 2(A \cup_1 A)A + 2A(A \cup_1A)$ and $$d(2A \cup_1(A\cup_1 A) = -4 (A\cup_1 A)\cup_1(A \cup_1 A) + 2A(A\cup_1A) - 2(A\cup_1A)A.$$
Adding shows that $$d((1/8)[A^2 + 2A\cup_1(A\cup_1 A]) = (1/2)[(A\cup_1 A)\cup_1 (A\cup_1A) + A(A\cup_1 A] .$$}
Thus the only essential cocycle change in $(1/2)x(a)$ would be to add  $(1/2)aSq^1a$.  This would result in the same change to $\D(w,p,a)$ as would result from changing the sign in the term $(\pm 1/4)A(A \cup_1A)$.  So we can fix $x(a).$  In \S4.5  just below we will clarify exactly how we want to fix $x(a)$.\\

It is clear that if the triple $(w,p,a)$ satisfies the basic equation for one choice of the sign in front of $(1/4)A(A \cup_1 A)$ in the basic equation, then the triple $(-w, p, a)$ satisfies the basic equation with the opposite sign.  So we fix this sign in the basic equation and the $k$-invariant, and work with the stated triples. Our choice of $x(a)$ and our choice of sign is motivated by the desire to relate the basic equations and the products  in $G^4(\Sigma X)$ and $G^3(X)$, under a suspension isomorphism  to be discussed in \S8.2.\\
 
\subsection{Anibal Medina's Formula and Our Formula for $x(a)$}

Medina's formula, communicated in a letter to us,  is $$x_M(a)(012345)  = a^2(01235) a^2(02345) \in C^5(\Z/2).$$  This is a remarkably simple formula, which indeed satisfies $dx_M = Sq^2Sq^2 + Sq^3Sq^1$ when applied to degree 2 cocycles. But to obtain better behaved formulas we must modify Medina's formula. The property of $x \in C^5(\Z/2)$ that we want is that it is the loop, or desuspension, of a natural cochain $\hat{x} \in C^6(\Z/2)$ that realizes the Adem relation  $d\hat{x} = Sq^2Sq^2 + Sq^3Sq^1$ when applied to degree 3 cocycles.  It turns out this can be accomplished by adding $aSq^1a$ to Medina's formula.  Thus, instead of Medina's formula, we take $x(a) = x_M(a) + aSq^1a$.
$$x(a)(012345) = a^2(01235) a^2(02345) + a(a\cup_1a)(012345)$$ $$ = a^2(01235) a^2(02345) + a(012)\big[ a(235)a(345)+a(245)a(234) \big] \in C^5(\Z/2).$$

In the next several sections we will return to the technical details of the construction of the group $G^4(X)$.  We need to define a product of triples $(w,p,a)(v,q,b)$ that satisfies the basic equation made explicit in \S\S4.2-4.5.  The class $x = x(a)$ defined here plays an important role in the basic equation, and we need to study that class carefully.\\

\section{Some Properties of the Class $x(a)$}

\subsection{Non-linearity of $x(a)$}

To define and study a product of triples, and to understand the relevant relations between triples, it is important to understand certain delicate properties of $x(a)$.  That is the purpose of this subsection and the next.  For someone browsing this paper, it would be reasonable to temporarily skip down to \S6.1 below, where the product is introduced. \\

Set $ \Delta x(a,b) = x(a+b) - x(a) - x(b)$ and set
\begin{align*}
\delta x(a, b)  = \ &  (a \cup_1a) \cup_1 (b\cup_1 b) + (a \cup_1 b)\cup_1(a\cup_1 b)\\
  & +\ (a\cup_1 a +b\cup_1b)(a\cup_2b) + (a\cup_2b)(a\cup_1a + b\cup_1b) \\
   & +\ (a^2+b^2)\cup_2 (a\cup_1 b) + (a\cup_1b)\cup_2 (a^2 + b^2) \\
   &  +\ (a\cup_1b)\cup_2 d(a\cup_1b) + (a\cup_2b)d(a\cup_2b) +\ a^2 \cup_3 b^2. 
 \end{align*}  
The result we are after is that $\Delta x(a, b) + \delta x(a,b) = dy_4(a, b)$ for a natural $\Z/2$ cochain $y_4(a, b)$ that is determined up to coboundaries by means of some further conditions. The term $y_4(a, b)$ occurs in \S6.2 below in the formula for the product of triples in $G^4(X)$.\\

We first claim that $\Delta x(a, b) + \delta x(a,b)$ is a cocycle.  Applying the formula $dx(a) = a^2 \cup_2 a^2 + (a \cup_1 a)^2$ to $a, b$ and $a+b$ and using the coboundary formula for cup$_i$ products, one can calculate $d\Delta x(a, b)$.  One can similarly compute $d\delta x(a, b).$  The computations agree, so $\Delta x(a, b) + \delta x(a,b)$ is a natural cocycle function of $a, b$, which must be cohomologous to a linear combination of $$aSq^1(a), bSq^1(b), aSq^1(b), Sq^1(a)b.$$  For this last statement, we appeal to the universal example $K(\Z,2 ,2) \times K(\Z/2, 2)$. But both $\Delta x(a,b)$ and $\delta x(a,b)$ vanish  when $a=0$ or $b=0$.  Clearly $\Delta x(a,a)$ vanishes and $\delta x(a,a) =a^2 \cup_3 a^2 = Sq^1(a^2)$, which is a coboundary. In fact, at the cocycle level it holds that $$Sq^1(a^2) = Sq^1(a)a + aSq^1(a) = (a\cup_1a)a + a(a\cup_1a) = d(a\cup_1(a\cup_1a)).$$ Thus the cocycle $\Delta x(a,b) + \delta x(a,b)$ must be a coboundary, or possibly a coboundary plus $Sq^1(a) b + aSq^1 (b) = Sq^1(ab)$, which becomes a coboundary in $\R/\Z$ cohomology since $Sq^1(ab)$ is the reduction of a torsion integral class. But we can resolve this last issue even in $\Z/2$ cohomology by brute force computation as follows.  All cup$_i$ products in the $\delta x(a, b)$ expression, evaluated on a 5-simplex, can be expanded as a sum of products of evaluations of the cocycles $a, b$ on faces of the simplex.  The $\Delta x(a, b)$ terms are directly given as such a sum, by Medina's modified formula in \S4.5.  We then wrote a computer program which found a $\Z/2$ 4-cochain $y_4(a, b)$ of this same form, a sum of products of evaluations of $a, b$ on faces of a 4-simplex, so that $dy_4(a, b) = \Delta x(a, b) + \delta x(a,b)$. \\ 

It is perhaps appropriate to include a few more words about our somewhat mysterious choice of $x(a)$ in \S4.5 and our ability to find a $\Z/2$ class $y_4(a, b)$ with $dy_4(a, b) = \Delta x(a, b) + \delta x(a,b)$.  In our early work on $G^4(X)$ we also developed some partial theory of the Pontrjagin dual of five dimensional Spin bordism so that there would be a suspension isomorphism $G^4(X) \to G^5(\Sigma X)$.  From that theory, it was clear that exactly one of the two classes, $x(a) = x_M(a)$ or  $x(a) = x_M(a) + a(a \cup_1 a)$ would be the desuspension of a natural 6-cochain $\hat x(\alpha)$ that satisfies the  Adem relation $d\hat{x}(\alpha) = Sq^2Sq^2(\alpha) + Sq^3Sq^1(\alpha)$ for $\Z/2$ 3-cocycles $\alpha$.  Moreover, for that $x(a)$, and not the other, there would  exist a natural $\Z/2$ class $y_4(a, b)$ with $dy_4(a, b) = \Delta x(a, b) + \delta x(a,b)$.  Specifically, that $y_4(a, b)$ would be the desuspension of a $\Z/2$ natural 5-cochain $y_5(\alpha, \beta)$ that satisfied an analogous formula $dy_5(\alpha, \beta) = \Delta \hat{x}(\alpha, \beta) + \delta \hat{x}(\alpha, \beta)$.  Computer computations within the $5D$ theory revealed that the class $x(a) = x_M(a) + a(a \cup_1 a)$ chosen in \S4.5 is the preferred class.  But this interesting $5D$ theory was not really needed in the end, since, as we stated above, we could just start with that  choice of $x(a)$  and find by computer  a $\Z/2$ class $y_4(a, b)$ with $dy_4(a, b) = \Delta x(a, b) + \delta x(a,b)$.  \\

We next want to state further conditions that pin $y_4(a,b)$ down up to coboundaries.  Since it holds that $dy_4(a,b) = 0$ when $a=0$ or $b=0$, and $dy_4(a, a) = Sq^1(a^2) $, we know $y_4(a, 0)$, $y_4(0, b)$, and $y_4(a, a) + a\cup_1(a\cup_1a)$ are cocycles.  By adding a linear combination of the cocycles $a^2, b^2, ab$ to $y_4(a, b)$, we could assume that $y_4(a, 0)$, $y_4(0, b)$, and $y_4(a, a) + a\cup_1(a\cup_1a)$ are coboundaries.  The last condition seems natural, but it turns out that in order for the product we will choose on $G^4(X)$ to be compatible with the evaluation formulas of \S3, we will need to add cocycles to $y_4(a, b)$ so that $y_4(a, a) + a\cup_1(a\cup_1a) +a^2$ is a coboundary. So we will assume $y_4(a, b)$ is chosen so that all three of the cocycles $$y_4(a,0),\ y_4(0, b),\  y_4(a, a) + a\cup_1(a\cup_1a) +a^2$$ are coboundaries. These three  coboundary conditions, along with the formula for $dy_4(a, b)$, imply $y_4(a, b )$ is well defined up to natural coboundaries.  In an appendix we will give an explicit formula for $y_4(a, b)$.\\

We repeat that the term $y_4(a, b)$ will occur  in the formula for the product of triples in $G^4(X)$.  Another product formula for triples is obtained by replacing $y_4(a, b)$ by $y_4(a, b) + ab$.  This is the same thing as requiring that $y_4(a, a) + a \cup_1 (a \cup_1 a)$ be a coboundary, rather than $y_4(a, a) + a\cup_1(a\cup_1a) +a^2$. As mentioned, the choice of $y_4(a, b)$ that we make here is the one that gives a product formula for triples that is compatible with the evaluation formulas of \S3.\\

\subsection{Value of  $x(a)$ on Coboundaries}

In order to determine the relations that will hold between triples in $G^4(X)$, we need to know the value of $x(a)$ on coboundaries.  Recall $x(dr) = x_M(dr) + dr(dr \cup_1dr)$, where $x_M(a)$ is Anibal Medina's cochain from \S4.5.  We claim
$$x_M(dr) = Sq^2(rdr) + Sq^1(r) Sq^1(dr) + dz_M(r).$$
where $z_M(r)$ is a natural cochain well defined up to exact cocycles.   The cochain operation $Sq^i(c) = c \cup_j c+ c \cup_{j+1}dc$, where $i + j = deg(c)$, commutes with the coboundary $d$. The claim then follows from the naturality argument since easily the two sides (without the $dz_M(r)$ term on the right) have the same image under $d$, so their difference $\Delta_M(r)$ is a 5-cocycle depending only on a $1$-cochain $r$.  It must thus be a coboundary. In fact, if $CX$ denotes the cone on $X$ then given $r\in C^1(X)$ we extend $r$ to $\hat r\in C^1(CX)$ by setting $\hat r$ to evaluate zero on every $1$-simplex in the cone $CX$ that is not contained in $X$.  Then we have the cocycle $\Delta_M(\hat r)$ on $CX$. This produces a suitable element $z_M(r)\in C^4(X)$ by $\langle z_M(r),\sigma^4\rangle=\langle \Delta_M(\hat r),C(\sigma^4)\rangle$.\\

Using the formula just above for $x_M(dr)$, the algorithm described in the paragraph above leads to an explicit  formula for $z_M(r)$ as a sum of products of evaluations of $r$ on edges of a 4-simplex.  We carried out this evaluation and the result is a formula for $z_M(r)$ that we will give in an appendix.\\

Once we find $z_M(r)$, we obviously have $$x(dr) =  Sq^2(rdr) + Sq^1(r) Sq^1(dr) + dz(r),$$ where $z(r) = z_M(r) + r(dr \cup_1 dr)$.  The term $z(r)$ occurs in the description of relations between triples in $G^4(X)$.

\section{The Product of Triples and a Non-Abelian Group}

\subsection{The Product of Triples}

We fix a simplicial complex $X$. Recall the definition from \S4.1 $$\widehat{\C} = \{(w, p, a) \in  C^4(\RR/\ZZ) \times C^3(\ZZ/2) \times Z^2(\ZZ/2)\ \vert \ dp = a^2\}.  $$ We define a product of triples in $\widehat{\C}$ by the formula $$ (w,p,a) (v,q,b) = (w+v+u,\ p + q +a\cup_1b,\ a+b)$$ 

\noindent where
$$u =  (1/2)p\cup_2q + (1/2) a^2 \cup_3 q + (1/2)(p+q)\cup_2(a\cup_1b) $$
$$  + (1/8)AB  - (1/4)(A\cup_1A) \cup_1 B - (1/4)(A + B)(A \cup_2 B) +(1/2)y_4(a,b).$$  Here, $y_4(a, b)$ is the $\Z /2$ cochain well-defined up to coboundaries discussed in \S5.1 above.  $A$ and $B$ are the special $\Z$ lifts of the $\Z/2$ cocycles $a$ and $b$.\\

\subsection{Additivity of the Basic Equation for the Product of Triples}

We claim that
$$ \D((w,p,a)(v,q,b)) = \D(w,p,a) + \D(v,q,b).$$

To prove this, recall the general definition from \S4.2, $$\D(w,p,a) = dw - [(1/2)Sq^2p + (1/4)A (A \cup_1A) + (1/2)x(a)],$$ where here $Sq^2p = p \cup_1 p + p \cup_2 a^2$ since $dp = a^2$.\\

The special lift  of the degree 2 cocycle $c = a+b$ is $C = A + B - 2(A\cup_2B)$, with $$dC = 2(C\cup_1C) =  2(A\cup_1A + B \cup_1 B - d(A\cup_2 B)).$$  It is then a lengthy but routine computation to see that a product  for triples of the form given in \S6.1 will be additive for $D$ exactly if $du=$
$$(1/2)[p\cup_1q + q\cup_1p +(p+q)\cup_1(a\cup_1b) + (a\cup_1b)\cup_1(p+q) +(1/2)[(p+q)\cup_2 d(a\cup_1b)]$$
$$ +(a\cup_1b)\cup_1(a\cup_1b) +(a\cup_1b)\cup_2(a^2+b^2)) + (a\cup_1b)\cup_2d(a\cup_1b)]$$
$$ +(1/4)[A(B\cup_1B) + B(A\cup_1A) - (A+B)d(A\cup_2B)]$$
$$+ (1/2)[ (a\cup_2b)(a\cup_1a)+ b\cup_1b + (a\cup_2b)d(a\cup_2b)]$$
$$+(1/2)[x(a+b) + x(a) + x(b)].$$
A direct computation shows that the formula for $u$ given in \S6.1 does indeed satisfy this equation.  The computation uses the coboundary formula for cup$_i$ products and the formula $dy_4(a, b) = \Delta x(a, b) + \delta x(a, b)$ from \S5.1.\\

We remark that since the Postnikov tower $E$ classifying $G^4$  is a loop space, one knows that a solution exists for a $u$ giving a product additive for $D$.  Then the above argument, when examined in detail, actually proves $(1/2)(\Delta x(a, b) + \delta x(a, b)) $ is a coboundary in $\RR/\ZZ$ cohomology.   The more delicate $\ZZ/2$ property, that  $\Delta x(a, b) + \delta x(a, b)$ is itself  the coboundary of some $\Z /2$ class $y_4(a, b)$, requires more discussion, as in \S5.1.\\

The product formula for triples is not explicit until an explicit formula for the cochain $y_4(a, b)$ is given.  Using a computer we have found a formula for $y_4(a, b)$ as a sum of 174 terms, each of which is a product of values of the 2-cocycles $a$ and $b$ on faces of 4-simplices.  For completeness, we write out an explicit formula for $y_4(a, b)$ in an appendix.\\
 
\subsection{A Non-Abelian Group}

With the product defined on $\widehat{\C}$ in \S6.1, it is clear that elements $(df, 0, 0)$ commute with all elements $(w, p, a) \in \widehat{\C}$.  Set $\bar{C}^4(\R/\Z) = C^4(\R/\Z) / dC^3(\R/\Z)$ and define  $$ \C = \widehat{\C}/{\{(df,0,0)\}} =  \{(w, p, a) \in \bar{C}^4(\R/ \Z) \times C^3(\ZZ/2) \times Z^2(\ZZ/2)\ | \  dp = a^2 \} .$$   We claim that $\C$ with the induced product  is a group, that is, associativity  holds and inverses  exist. In fact, inverses are given by $$(w,p,a)^{-1} = (w', p+Sq^1a, a)$$ where $$w' = -w -(3/8)A^2 -(1/4)(A\cup_1A)\cup_1A +(1/2)Sq^1a\cup_2p .$$ Verification that this gives the inverse relies on the rather subtle fact that $(1/4)A^2 + (1/2)Sq^1p$ is a coboundary, and also that $y_4(a, a) + (a\cup_1a)\cup_1a + a^2$ is a coboundary.  Verification that associativity holds in $\C$ follows the lines of the proof of a corresponding result in the $3D$ theory in ([1], \S3.4).\\

The map $\D$ of \S4.2 and \S6.2 descends to a function $\D\colon \C \to C^5( \RR / \ZZ)$, which is  a group homomorphism. Triples $(w,p,a) \in Kernel(\D)$ name simplicial maps $X \to E$, where $E$ is the Postnikov tower of \S4.3 that represents the dual of reduced $4D$ Spin bordism.  In \S7 we will study the  triples  $(w,p,a) \in Kernel(\D)$ that represent null-homotopic maps.  We can formulate these null-homotopic triples as the image of another homomorphism $\D' \colon \C' \to \C$, so that $Image(\D')$ contains all commutators.  Then we define the main object of the paper, the abelian group $G^4 = Kernel(\D) / Image (\D')$.\\

\subsection{Variations of the Product Formula} 

At this point it is appropriate to point out that there are alternate possible formulas for the product of triples in \S6.1 above that lead to an associative group $\C$ and an abelian group $G^4 = Kernel(\D) / Image (\D')$. Note that additivity of the basic equation in \S6.2  forces $du$.  It certainly changes nothing to add a coboundary to $u$, since for triples in $\C$ the first coordinate is only defined up to coboundaries.  Thus the cochain product formula for triples is extremely non-unique.  But one could also try to add to $u$ a linear combination of cocycles $(1/4)\mathcal{P}(a), (1/4)\mathcal{P}(b), (1/2)ab$, where $\mathcal{P}$ is the Pontrjagin square.  Only adding $(1/2)ab$ is consistent with $(0,0,0)$ representing the identity element.  Adding $(1/2)ab$ does give another associative product. But the result of adding $(1/2)ab$ to $u$ yields a group isomorphic to $\C$, if the triples $(w,p,a)$ are corresponded to the triples $(w+(1/4)\mathcal{P}(a), p, a)$.  The fact that we have characterized $y_4(a, b)$ up to coboundaries in \S5.1 thus pins our product down up to adding coboundaries, and possibly adding $(1/2)ab$.\\

Notice that adding $(1/2)ab$ to $u$ is the same as just adding $ab$ to $y_4(a, b)$ and leaving the rest of $u$ alone.  Naming one of these two choices of  $y_4(a, b)$, up to coboundary,  can also be expressed as asserting which of the cocycles $y_4(a, a) + a \cup_1 (a \cup_1 a)$ or $y_4(a, a) + a \cup_1 (a \cup_1 a) + a^2$ is a coboundary.  In \S5.1 we made the latter choice, so therefore with that choice of $y_4(a, b)$ we have explicitly singled out one of the two possible products on $\C$.

\section{Relations Between Triples and  the Group $G^4(X)$}
 
\subsection{Relations Between Triples}

First Relations:  $(df, 0, 0) \equiv (0,0,0)$  where $f \in C^3(\RR/\ZZ)$\\

\noindent Second Relations: $((1/2)Sq^2c,\ dc,\ 0) \equiv (0,0,0)$  where $c \in C^2(\ZZ/2)$ \\

\noindent Third Relations:  $((1/2)z(r)+(1/4)R (Dr \cup_1 Dr),\ r dr,\ dr) \equiv (0,0,0)$ where $r \in C^1(\ZZ/2)$\\

(Here $R$ is the special lift of $r$ and $Dr$ is the special lift of $dr$. Note $dDr = 2 Dr \cup_1 Dr $. The term $z(r)$ is the term discussed in \S5.2 above.  We give an explicit formula for $z(r)$ in an appendix.)\\

The relations are all in $Kernel(\D) \subset \C$. Verification that the third relation is in $Kernel(\D)$ makes use of  the formula for $dz(r)$ in \S5.2, along with the formula $(1/4)dR =  (1/4)Dr + (1/2)Sq^1(r) $, which is proved by direct evaluation on 2-simplices.\\

There is a systematic method for finding relations.  Triples $(w,p,a)$ satisfying the basic equations $da = 0, dp = a^2$ and $dw = k(p,a)$ are regarded as simplicial maps $X\to E$, where $E$ is the Postnikov tower defined in \S4.3.  A triple represents a null-homotopic map exactly when there is an admissible triple $(\hat{w}, \hat{p}, \hat{a})$ of cochains on the cone $CX$ that restricts to $(w,p,a)$ on the base $X$. By exploiting a chain level null-homotopy of $X$ in $CX$, one can find formulas for the restrictions of $\hat{w}, \hat{p}, \hat{a}$ to the base.  These formulas give the relations.\\

Our construction of $G^4(X)$ uses cochain triples on a fixed ordered simplicial structure on $X$.  With that goal, we need to write down explicit relations, as above, in terms of a single ordered simplicial structure on $X$.  But there is another important view of relations.   The group $G^4(X)$ we construct  is a homotopy functor.  So if $I \times X$ is given an ordered simplicial structure, then the inclusions of the ends induce  isomorphisms $G^4(I \times X) \simeq G^4(\{i\} \times X)$, for $i = 0, 1$. The two ordered simplicial structures on $\{0\} \times X$ and $\{1\} \times X$ can be different.  But it now makes sense to say triples $(w_0, p_0, a_0) \in G^4(\{0\} \times X)$ and $(w_1, p_1, a_1) \in G^4(\{1\} \times X)$ represent the same element of $G^4(X)$.  Specifically, this will be the case if there is an element $(\hat{w}, \hat{p}, \hat{a}) \in G^4(I \times X)$ that restricts to $(w_i, p_i, a_i)$ on the ends.  This point of view will be quite important in \S9, where we discuss computer methods used to make a crucial computation in $G^4( S^2 \times S^2)$.  In theory this computation could have been carried out with relations in the $G^4$ group of the second barycentric subdivision of a standard triangulation of $S^2 \times S^2$.  But those computations were too big.  Instead we were able to work $I \times (S^2 \times S^2)'$, where $(S^2 \times S^2)'$ is the first barycentric subdivision, with different (ordered!) simplicial structures on the two ends.\\

\subsection{Amalgamation of Relations} 

Set $\C' = C^2(\ZZ/2) \times C^1(\ZZ/2)$. Define $\D' \colon \C' \to \C$ by $\D'(c, r) =$
$$ ((1/2)dc\ \cup_2\ rdr + (1/2)Sq^2c +(1/2)z(r)+ (1/8)R d(Dr),\ dc + rdr,\ dr).$$
Then  $\D \circ \D' = 0$.\\

\subsection{Products of Relations} 

Define a product on $\C'$ by $(c,r)(e,s) = (c + e + r\cup_1 ds + sr,\ r + s)$.  Then $$\D'((c,r)(e,s)) = \D'(c,r)\D'(e,s) \in \C.$$   $\C'$ is a group and $Image(\D') \subset \C$ is a subgroup containing all commutators. In fact, $\{\D'(c,0)\}$  contains all commutators. Proofs of the product formula for relations and the statements about commutators follow the lines of corresponding proofs in ([1], \S3.3).\\

\subsection{Definition of the Group $G^4(X)$}

We now have for any simplicial complex or simplicial set $X$ homomorphisms of  groups with $\D \circ \D' = 0$, 
$$ \C'(X) \xrightarrow{\D'} \C(X) \xrightarrow{\D}  C^5(X, \R/\Z).$$
Moreover $Image(\D')$ contains all commutators in $\C$. Therefore we can define an abelian group, functorial in  $X$,
$$G^4(X) = Kernel(\D) / Image (\D')$$
In \S3 we have outlined the construction of an evaluation map from $G^4(X)$ to the Pontrjagin dual of $\widetilde{\Omega}^{spin}_4(X)$. The first two relations are seen to evaluate as 0.  The third relation is different, and the evaluation map simply declares these relations to evaluate as 0. The proof that the construction is additive in the $G^4$ variable is rather subtle and indirect, as explained in \S3. After that, the proof that the evaluation is an isomorphism $$G^4(X) \simeq {\rm Hom}(\widetilde{\Omega}^{spin}_4(X), \R / \Z)$$ follows from an argument comparing natural  filtrations on both sides.  These filtrations are described in the next section.\\

\subsection{The Subgroups and the Filtration Quotients} 

Set $SSH^2(X; \ZZ/2) = \{\bar{a} \in H^2(X; \ZZ/2)\ \vert\ \bar{a}^2 = 0, \  < (1/2)Sq^2, Sq^2 > \bar{a}  = 0\}$, where $<\cdot, \cdot>$ is the indicated secondary operation.\\

\noindent Set $SH^3(X; \Z/2) = \{\bar{p} \in H^3(X; \ZZ/2)\ \vert \ (1/2)Sq^2(\bar{p}) = 0\}$.\\

\noindent Set $QH^4(X; \RR/\ZZ) = H^4(X; \RR/\ZZ) / \{(1/2)Sq^2(H^2(\ZZ/2))\}$.\\

There is a natural filtration on $G(X) = G^4(X)$:
$$G(X)=G^0(X)\supset G^1(X)\supset G^2(X)\supset 0$$
where $G^1(X)$ and $G^2(X)$ are the  subgroups consisting of all elements of $G(X)$ represented by triples of the form $(w,p,0)$ and $(w,0,0)$, respectively. The following results are proved in the same manner as the $3D$ analogues in our previous work in ([1], \S4.1).
\begin{enumerate}
\item the map $(w,p,a)\mapsto a$ determines a natural  isomorphism of  abelian groups
$$G(X)/G^1(X)\to SSH^2(X;\Zee/2\Zee),$$ 
\item the map $(w,p,0)\mapsto p$ determines a natural isomorphism of abelian groups 
$$G^1(X)/G^2(X)\to SH^3(X;\Zee/2\Zee),$$

\item  the map $(w,0,0)\mapsto w$ determines a natural isomorphism of  abelian groups
$$G^2(X)\to QH^4(X;\Ar/\Zee).$$
\end{enumerate}
These filtration quotients are identical to the filtration quotients that arise from filtering the Pontrjagin dual of reduced $4D$ Spin bordism of $X$ by $\R / \Z$ valued homomorphisms vanishing on  images of reduced Spin bordism of skeletons  $X^{(2)} \subset X^{(3)}  \subset X$. This claim follows immediately from knowledge of the differentials in the Atiyah-Hirzebruch spectral sequence for the Pontrjagin dual of Spin bordism.\\

But there is a very important point to make about these cohomology groups.  In the context of the Atiyah-Hirzebruch spectral sequence, the indicated cohomology groups are really the Pontrjagin duals of homology groups $H_i(X; \widetilde{\Omega}^{spin}_j(pt))$, where $i+j = 4$.  The subgroups indicated by $S$ and $SS$ are duals of quotients of homology groups.  The quotient indicated by $Q$ is dual to a subgroup of a homology group.\\
  
The group called $G^1(X)$ here is the group called $\tilde{G}(X)$ in \S3.4 and \S10.3.  In \S10.3, we give the filtration argument that proves that evaluation $$G^1(X) \to  {\rm Hom}(\widetilde{\Omega}^{spin}_4(X) / Image(\widetilde{\Omega}^{spin}_4(X^2)), \R / \Z)$$ is an  isomorphism.  This is essentially a comparison of exact sequence $(I)$ below with the dual of a corresponding exact sequence arising from the Atiyah-Hirzebruch bordism spectral sequence filtration.  To finish the proof that the evaluation $$G^4(X) \to {\rm Hom}(\widetilde{\Omega}^{spin}_4(X), \R / \Z)$$ is an isomorphism, we need to make a similar comparison of exact sequence $(III)$ below with the dual of an Atiyah-Hirzebruch  bordism exact sequence associated to the filtration $$ Image(\widetilde{\Omega}^{spin}_4(X^1)) \subset Image(\widetilde{\Omega}^{spin}_4(X^2))  \subset \widetilde{\Omega}^{spin}_4(X).$$   Note  $\widetilde{\Omega}^{spin}_4(X^1) = 0$ since reduced 3-dimensional Spin bordism of a point is 0.   \\

Now, the $G^1(X)$ part is already done.  The last step then is to regard a Spin bordism class $M^4 \to X^2 / X^1$ as a map from $M^4$ to a wedge of 2-spheres.   In \S3.3, we observed that the evaluation of a map $g \colon M^4 \to S^2$ on the generator $(0,0, s\alpha)$ of $G^4(S^2)$ is the Arf invariant in $\Z/2$ of a framed surface in $M$ obtained as $g^{-1}(pt)$.  The Arf invariant is exactly the Spin bordism class of the framed surface.  This observation completes the proof that exact sequence $(III)$ below maps isomorphically to the dual of the corresponding Atiyah-Hirzebruch filtration exact sequence, and hence completes the proof that our evaluation $$G^4(X) \to {\rm Hom}(\widetilde{\Omega}^{spin}_4(X), \R / \Z)$$ from \S3 is an isomorphism.\\  

As is well-known, the $E_\infty$ page of the Atiyah-Hirzebruch spectral sequence determines filtration quotients, but does not  determine the group extensions below.   
$$(I) \hspace{.5in}0 \to QH^4(X;\Ar/\Zee) \to G^1(X) \to SH^3(X;\Zee/2\Zee) \to 0 \hspace{.5in}$$
$$ (II)\hspace{.2in} 0 \to SH^3(X;\Zee/2\Zee)  \to G(X)/ G^2(X) \to SSH^2(X;\Zee/2\Zee) \to 0 \hspace{.2in}$$
$$ (III) \hspace{.7in} 0 \to G^1(X) \to G(X) \to SSH^2(X;\Zee/2\Zee) \to 0 \hspace{,7in}$$
$$ (IV)\hspace{.5in}  0 \to QH^4(X;\Ar/\Zee) \to G(X) \to G(X) / G^2(X) \to 0 \hspace{.6in}$$
Since we work at the cochain and cocycle level to describe all the middle groups in the four sequences, we can write down formulas that determine these group extensions up to isomorphism. In Appendix \S10.2 we will give further details about these extensions.

\section {Applications of Suspension of Cochains} 

\subsection{Suspension of Cochains}

We regard the suspension $\Sigma X$ of a space to be the obvious  union of two cones $C^{+}X$ and $C^{-}X$.  Given a triangulation of $X$ with vertex order,  we label the new upper cone vertex $+ \infty$ and the new lower cone vertex $-\infty$.  That is, the cone vertices always occur as the last or first vertex of simplices in the suspension.  But we won't use the lower cone vertex. Given a cochain $c \in C^n(X)$, with any coefficients, we define $s(c) \in C^{n+1}(\Sigma X)$ as follows.  On any simplex in the lower cone $C^{-}X$, the value of $s(c)$ will be 0.  On a simplex of form $(012 ... n \infty)$ in the upper cone $C^{+}X$, the value will be $s(c)(01 ... n \infty) = c(01 ... n)$.\\

The `cone vertex last' convention results in the easily proved formula $sd = ds$.  That is, $s$ is a cochain map $s \colon C^*(X) \to C^{*+1}(\Sigma X, C^- X)$.   The suspension cochain map  $s$ induces suspension isomorphisms on cohomology with any coefficients $\tilde H^*(X) \simeq H^{*+1}(\Sigma X, C^-X) = H^{*+1}(\Sigma X)$.\\

The suspension $s$ has some very nice properties relating $\cup_i$ products in $X$ and $\Sigma X$. First, we point out that with the given ordered triangulation of $\Sigma X$, all ordinary cup products $sx \cup_0 sy$ are 0.  The reason is, an ordered simplex can have at most one vertex $+\infty$, so a proper `first face' will always lie in $C^{-}X$.  On the other hand, the following remarkable formula holds for all $i \geq 0$:  $$s (x \cup_i y)= (-1)^{deg(x)+i+1}sx \cup_{i+1} sy.$$
We believe this is an important formula.  It is not easy to prove.   Obviously it implies  that Steenrod square operations commute with suspension, not just on cohomology and cocycles, but actually on all cochains.  The cochain cup$_i$ operations generalize to other multi-variable cochain operations and there should be useful extensions of this suspension formula to these other operations.\footnote{In hindsight, we believe using the lower cone $C^-(X)$ to define a suspension $s$ with $ds = -sd$ is the more natural choice. But then in order to get the cleanest formulas relating $s$ and cup$_i$ products it is necessary to use alternatives to the historical definitions of cup$_i$ products, including ordinary cup product.  But it is hard to overturn historical conventions!}\\

\subsection{The Suspension $G^3(X) \to G^4(\Sigma X)$}

We define $s \colon G^3(X) \to G^4(\Sigma X)$ by $s(w,p,a) = (sw, sp, sa).$  The first thing to check is that if $(w,p,a)$ satisfies the $G^3$ basic equations $dp = 0$ and  $dw = (1/2)p^2$ then $(sw, sp, sa)$ satisfies the $G^4$ basic equations $dsp = (sa)^2$ and    $D(sw, sp, sa) = 0$.  But $sp$ is a cocycle and ordinary cup products of suspension classes vanish, so $(sa)^2 = 0$.   Also since products of suspensions vanish, from \S4.5 we see $x(sa) = 0$.  For the same reason, the term  $(1/4)sA(sA \cup_1 sA)$ in the $G^4$ basic equation also vanishes.  This leaves only $$dsw = sdw = s((1/2)p^2) = (1/2)s(Sq^2p) = (1/2)Sq^2(sp),$$ as desired.\\

Next, in $G^4$ many of the terms of a product of suspension triples vanish, leaving only $$(sw,sp,sa)(sv, sq, sb) = (sw+sv+\hat{u}, sp+sq+sa\cup_1sb, sa+sb)$$ where $$ \hat{u} = (1/2)(sp\cup_2sq) +(1/2)(sp+sq)\cup_2(sa\cup_1 sb)$$ $$-(1/4)(sA \cup_1 sA) \cup_1 sB +(1/2)y_4(sa, sb).$$

Notice the term $(-1/4)(sA \cup_1 sA) \cup_1 sB$,  which will help tell us how the $4D$ multiplication desuspends to $3D$. The sign of that term is delicate. Since $A$ has degree 1 here, the suspension formula in \S8.1 gives $$s(A^2B) = (-1) s(A^2) \cup_1 sB = (-1)(sA \cup_1 sA) \cup_1 sB.$$ One then sees $$\hat{u} = (1/2)s(p\cup_1q) + (1/2)s((p+q)\cup_1ab) + (1/4)s(A^2B) + (1/2)y_4(sa, sb).$$
We must still evaluate $(1/2)y_4(sa, sb)$. From \S5.1, and the fact that $ x(sa) = x(sb) = x(sa + sb) = 0$, we see that most terms in $dy_4(sa, sb)$ are, or contain, products of suspensions, leaving only $$dy_4(sa, sb) = s(a^2 b^2) + s(abab).$$
But $a^2 b^2 + abab = d(a(a\cup_1 b)b)$, so we have $dy_4(sa, sb) = ds(a(a\cup_1b)b)$.  Therefore, modulo cocycles, we have  $$y_4(sa, sb) \equiv s(a(a\cup_1b)b).$$
The left side is a coboundary  when $a = 0$ or $b = 0$, and the right side vanishes. When $a = b$, on the left we know from \S5.1 that $$y_4(sa, sa)+ sa\cup_1(sa\cup_1sa) + (sa)^2 = y_4(sa, sa) + s(a^3)$$ is a coboundary.  On the right, when $a=b$ we  get $s(a(a\cup_1b)b) = s(a^3)$.   The conclusion is that up to coboundaries  $$y_4(sa, sb) = s(a(a\cup_1b)b ).$$  
 Backing up, we have shown   $$\hat{u} = (1/2)s(p\cup_1q)) + (1/2)s((p+q)\cup_1ab) + (1/2)s(a(a\cup_1b)b) + (1/4)s(A^2B).$$  So the multiplication in $G^4$ desuspends to the multiplication in $G^3$  as claimed in \S2.2.  In fact, the argument shows that both products on $G^4$ discussed in \S6.4 desuspend to the multiplication in $G^3$.\\
 
\subsection{The Factorization of Elements of $G^4(M)$}
 
 The first few homotopy groups of $\Sigma \R P^\infty$ are $\pi_2 = \Z/2,\ \pi_3 = \Z /2$ and $ \pi_4 = \Z /4$.  The first $k$-invariant `kills' $(\bar{s\alpha})^2$, where $\alpha \in Z^1(\R P^\infty; \Z/2)$ is a cocycle generating $H^1(\R P^\infty; \Z/2)$, and $s$ is the suspension of \S8.1. It is also possible to write down the second $k$-invariant, but we won't need this.  We fix an ordered triangulation of $\R P^\infty$ and then the double cone geometric model of $\Sigma \R P^\infty$, as in \S8.1.   Since cup products of suspensions vanish at the cochain level in this model, it follows that the triple $(0, 0, s\alpha)$ represents an element of $G^4(\Sigma \R P^\infty)$.  That is, $\D(0, 0, s\alpha) = 0$. Specifically, referring to the definition of $\D$ in \S4.2 along with the definition of $x(a)$ in \S4.5,  both the special lift term $(1/4)S\alpha (S\alpha \cup_1 S\alpha) $ and the term $(1/2)x(s\alpha)$ are zero.\\
 
A simple first conclusion is that given a Spin manifold $M^4$ and a class $\bar{a} \in H^2(M^4; \Z/2)$, elementary obstruction theory and the fact that  $\bar{a}^2 = 0$ implies there are maps $u\colon M^4 \to \Sigma \R P^\infty$ so that $\bar{s\alpha}$ pulls back to $\bar{a}$.  We can triangulate $M$ so that $u$ is an ordered simplicial map.  Then whatever the cocycle $a \in Z^2(M; \Z/2)$, there will be a 1-cochain r so that $a + dr = u^*(s\alpha).$ Hence $(0, 0, a+dr) \in G^4(M)$.\\

Now consider $(w, p, a) \in G^4(M)$.  Multiply by a third relation from \S7.1, specifically the relation $((1/2)z(r)+(1/8)Rd(Dr), rdr, dr)$,  resulting in an equivalent element $(w, p, a) \equiv (w', p', a')$, where $a' = a+dr$.  Then we also see $(w', p', a') = (w', p', 0)(0, 0, a') \in G^4(M)$.  This is the factorization of $(w, p, a)$ that is used in \S\S3.1-3.3 to define the evaluation $\langle (w, p, a), (M \xrightarrow{Id} M) \rangle$.\\

\section {The Evaluation $G^4(X) \to {\rm Hom}(\widetilde{\Omega}_4^{\rm spin}(X),\Ar/\Zee)$}

\subsection {The Uniqueness of the Evaluation } 
   
As we pointed out in \S6.4, there are two products on the underlying set $G^4(X)$, which define isomorphic groups, via the correspondence of triples given by $(w,p,a) \longleftrightarrow (w+ (1/4)\mathcal{P}(a), p, a)$, where $\mathcal{P}$ denotes the Pontrjagin square.  We have chosen one product, with our choice of $y_4(a, b)$, and we claim that the resulting group $G^4(X)$ admits no non-trivial natural automorphisms other than the inverse automorphism of abelian groups.  This means that  up to sign there will be a unique natural evaluation isomorphism $G^4(X) \to {\rm Hom}(\widetilde{\Omega}_4^{\rm spin}(X),\Ar/\Zee)$. The sign can be pinned down by agreeing that elements $(w,0,0)$ evaluate on $f \colon M^4 \to X$ as $\int_M f^*w$.  This is part of the definition of the evaluation pairing discussed in \S3.4.\\

We explain here why the automorphism group of the functor $G^4(X)$ consists only of the identity and the inverse map.  We work directly with triples $(w, p, a)$.  Suppose $\phi \colon G^4(X) \to G^4(X)$ is a natural automorphism, say $\phi(w,p,a) = (w', p', a')$.  Then obviously the cohomology classes $\bar{a}, \bar{a}'$ are the same, so $a'  = a + dx$.  Multiplying by a third relation from \S7.1, which does not change the automorphism, we can assume $a' = a$.  Next we must have either $p' = p + dc$ or $p' = p + dc + Sq^1a$, since up to coboundaries $Sq^1a$ is the only non-zero natural 3-cocycle function of a 2-cocycle $a$.  Multiplying by a second relation from $\S7.1$, we can get rid of the $dc$ term.  Then, by composing $\phi$ with the inverse automorphism from \S6.3 if necessary, which has the form $(w, p, a) \mapsto (*, p+Sq^1a, a)$, we can assume in completing the proof that $p' = p$.  Finally, we must have $dw = dw' = k(p, a)$, from the basic equation in \S6.1.  But there are no natural non-zero candidates for $$\bar{w-w'} \in H^4(K(\Z/2, 2) \ltimes_{a^2} K(\Z/2, 3)\ ;\ \R / \Z)$$ other than the order 2 element $(1/4)\mathcal{P}(a)$.  But $\phi(w, p, a) = (w + (1/4)\mathcal{P}(a), p, a)$ is not a group automorphism, because $$(1/4)\mathcal{P}(a+b) = (1/4)\mathcal{P}(a) + (1/4)\mathcal{P}(b) + (1/2)ab.$$  (We have already seen that this map $\phi$ transforms one of the products on the set $G^4(X)$ to the other product.) Thus $w'$ and $w$ differ by a coboundary, hence $(w', p', a') =  (w, p, a) \in G^4(X)$.\\

 With the sign choice two paragraphs above, we have  now  specified a unique natural evaluation homomorphism $G^4(X)$ $\to {\rm Hom}(\widetilde{\Omega}_4^{\rm spin}(X),\Ar/\Zee)$.  In \S3 we described necessary formulas for this evaluation, up to a single ambiguity in the evaluation of $(0,0, s\alpha)$ on a Spin bordism element $g\colon M^4 \to \Sigma \R P^\infty$.  We will settle this ambiguity by taking $M = S^2 \times S^2$.  A tubular neighborhood of the diagonal is the tangent disk bundle of $S^2$, with boundary $\R P^3$.  We thus have a map $g\colon S^2 \times S^2 \to \Sigma \R P^\infty$, transverse to $\R P^\infty$, which maps the tangent disk bundle to the lower cone $C^-(\R P^\infty) \subset \Sigma \R P^\infty$ and the complement of the tangent disk undle to the uppper cone $C^+(\Ar P^\infty)$.  From \S3.3, the orientation we choose on $N = \R P^3$ here is the boundary of the inverse image of the upper cone, which is the negative of the standard orientation.  Therefore, from the discussions in  \S2.2 and \S3.3, the result we want is that $$\langle (0, 0, s\alpha),\  [S^2 \times S^2 \xrightarrow{g} \Sigma \R P^\infty] \rangle = \rm{Arf}(-\R P^3) = -1/8.$$ In the next section we will study $G^4(S^2 \times S^2)$ further and prove that this evaluation must indeed be $-1/8$.\\

\subsection  {A Product Formula on $S^2 \times S^2$}
Here we continue to deal with ordered triangulations meaning triangulations whose vertices are partially ordered with the partial order restricting to a total ordering of the vertices of any simplex. Given an ordered triangulation we have
the complex of normalized simplicial cochains, those that vanish on degenerate ordered simplices.\\

Let $S^2 \times S^2$ have the natural product orientation from the usual orienta-
tion on $S^2$. This determines a Spin structure on $S^2 \times S^2$. We begin with the
product cell decomposition of $S^2 \times S^2$ coming from triangulations of each
factor as the boundary of the 3-simplex. We then take the triangulation of
$S^2 \times S^2$ whose vertices are the 0-cells of this cell decomposition, hence
each vertex is an ordered pair of integers, each integer being between 0 and
3 inclusive. We order these vertices lexacographically.  A set of vertices is the set of vertices of a simplex of this
triangulation if and only if (i): they all lie in the closure of the same product cell,
and (ii): under the lexacographic ordering of the vertices both the induced
ordered set of first integers and the induced ordered set of second integers is
weakly ordered in the usual ordering on integers. Notice that the diagonal
copy of $S^2$, denoted $D$, is a full subcomplex of this triangulation.
We let $T$ be the first barycentric subdivision with the vertex partial
ordering given by assigning to a barycenter of a simplex $\sigma$ the dimension of $\sigma$.
The union of all closed simplices of $T$ that meet $D$ form a regular
neighborhood $V$ of $D$ and its boundary $\partial V$ is PL homeomorprhic to $\R P^3$.  We identify $\partial V$ with $\R P^3$.\\

The projections, $\pi_1$ and $\pi_2$, of $S^2 \times S^2$ onto the two factors are order preserving simplicial maps from $T$ to the first barycentric subdivision of $\partial \Delta^3$. Let $U \in C^2(\partial \Delta^3; \Z)$ be the cochain that takes value 1 on a single non-degenerate 2-simplex and is zero on all other 2-simplices. It is a cocycle
whose cohomology class  generates $H^2(\partial \Delta^3; \Z)$. 
 Let $A_1 = \pi_1^*U$ and $A_2 = \pi_2^*U$ 
 be the pullback cocycles in $C^2(S^2 \times S^2; \Z)$. 
 These form a basis for $H^2(S^2 \times S^2; \Z)$.
 Furthermore, by naturality of the cup product, for i = 1, 2,
we have $A_i^2
 = 0$ and $A_i \cup_1 A_i = 0$. Of course, $A_1\cup A_2$ evaluates 1 on the
fundamental cycle of $S^2 \times S^2$, Let $a_i$ be the $\Z/2$ reduction of $A_i$. It follows
that the triple $(0, 0, a_i)$ represents an element of $G^4(S^2 \times S^2)$. Since the
preimage under $\pi_i$ of a point in $\partial \Delta^3$ is a two-sphere, this preimage is a Spin
boundary and hence the element of 
$\Omega^{spin}
_4 (\partial \Delta^3)$ represented by $(S^2 \times S^2, \pi_i)$
is 0. Again by naturality, it follows that the elements $(0, 0, a_i) \in G^4(S^2 \times S^2)$ evaluate 0 on the identity map $Id \colon S^2 \times S^2 \to S^2 \times S^2$.\\

By the product formula
$$(0,0,a_1)(0,0,a_2) = ((1/8)A_1A_2 + (1/2)y_4(a_1, a_2),\ a_1 \cup_1 a_2,\ a_1 +a_2)$$
Direct computer computations shows that our choice of $y_4(a_1, a_2)$ from  \S5.1 integrates to 0 over $S^2 \times S^2$.
Hence the first coordinate of $(0,0,a_1)(0,0,a_2)$ integrates to $1/8$ on the fundamental
class of $S^2 \times S^2$.  This is the first key computer result.\\

\subsection {The Diagonal Class}

In this subsection we deal with the same triangulation of $S^2 \times S^2$ as in
the previous subsection, but with a different vertex partial-ordering. Each
vertex in the interior of $V$, which are the vertices of the diagonal $D$, is assigned 5 more than the dimension of the
simplex of which it is the barycenter, each vertex in $\partial V = \R P^3$ is assigned 10 more
that the dimension of the simplex of which it is the barycenter, and each
vertex in the complement of $V$ is assigned 15 more than the dimension of
the simplex of which it is the barycenter. This means that each vertex in the
interior of $V$ is less then each vertex in $\partial V$, and the latter are less than each
vertex in the complement of $V$. This triangulation with its vertex partial
order is denoted $T'$. We are identifying $\partial V$ with $\R P^3$. There is an order-preserving
simplicial map $g\colon T' = S^2 \times S^2 \to\Sigma \R P^3 \subset \Sigma\R P^\infty$
sending $V$ to the lower cone and the complement of $V$ to the upper cone
with the property that the restriction to the middle level is the identity map $Id\colon \partial V \to \R P^3$.\\
 
Let $\alpha \in Z^1(\R P^\infty; \Z/2)$  be a cocycle whose cohomology class is non-trivial. Let
$s\alpha \in Z^2(\Sigma \R P^\infty; \Z/2)$ denote the 2-cocycle obtained as in \S8.1 by suspending $\alpha$
over the postive cone and being identically zero on the negative cone. This
cocycle generates the second cohomology of $\Sigma \R P^\infty$ and has square zero,
as a cocycle. 
The pullback $c = g^*(s\alpha)$ is a $\Z/2$-cocycle of square zero
representing the $\Z/2$ cohomology class $a_1+a_2 \in H^2(S^2 \times S^2; \Z/2)$.
It follows that $(0,0, c)$ represents an element
in $G^4(S^2 \times S^2)$.\\

In general, the functor $G^4(X)$ does not depend on the choice of ordered triangulation of $X$. Namely, one can work with ordered triangulations of $X \times I$ that extend two given  triangulations on $X \times \partial I$, and exploit homotopy invariance of the functor $G^4(X)$. For the triangulation $T'$ of $S^2 \times S^2$, there will be an element $(w', p', a') \in G^4(S^2 \times S^2)$ that represents the same element as the product $(0,0,a_1)(0,0,a_2)$ represents in the triangulation $T$.  For any such choice, $a' + c$ will be a coboundary, since both $a'$ and $c$ represent the homology class of $a_1 + a_2$.  Multiplying $(w', p', a')$ by a third relation from \S7.1 if necessary, we may assume $a' = c$.  Then $dp' = c^2 = 0$.  Since $H^3(S^2 \times S^2; \Z/2) = 0$, we know $p'$ is a coboundary and we can multiply by a second relation from \S7.1 if necessary and assume $p' = 0$. Thus for the triangulation $T'$, we can assume $(w',p',a') = (w, 0, c) = (w,0,0)(0,0,c)$.  The cocycle $w$ is well-defined up to coboundaries. \\ 

CLAIM: $(w,0,0)$ evaluates as $+1/8$ on the Spin bordism element of $S^2 \times S^2$
represented by the identity map, and $(0,0,c)$ evaluates as $-1/8$.\\

In order to prove this claim, we must find a way to compare $(w,0,c)$ and $(0,0,a_1)(0,0,a_2)$. We
cannot do this directly since the elements are defined with respect to different partially order triangulations. We make the comparison in the following way. We consider $S^2 \times S^2 \times I$. The vertices of the triangulation
of this product are the vertices of $T$ at the 0-end union the vertices of $T'$
at the 1-end, each with their function to the non-negative integers. We
triangulate the product with this set of vertices in the usual way. That is to say the vertices of the simplices of dimension 5 are those of the form $((v_0,0), \dots , (v_i, 0), (v_i, 1), \dots , (v_4,1))$ where the 
$(v_0, \dots, v_4)$ range over the ordered vertices of 4-simplices of $T$ and $0 \leq i \leq 4$. This produces a triangulation of the product and a partial ordering of the vertices extending $T \times \{0\}  \coprod T' \times \{1\}$. Denote the resulting ordered
triangulation by $\widetilde T$.\\

Computer computations allowed us to extend the $\Z/2$-cocycle on $S^2 \times S^2 \times \partial I$
given by $a_1 + a_2$ on $T \times \{0\}$ and by $c = g^*s\alpha$ on $T' \times \{1\}$ to a $\Z/2$-cocycle $\tilde a$ on $\widetilde T$ and then to
extend the 3-cochain on $S^2 \times S^2 \times \partial I$  given by $a_1 \cup_1 a_2$ on $T \times \{0\}$ and 0 on $T' \times \{1\}$ to a $\Z/2$ cochain $\tilde p$ such that $d\tilde p = \tilde{a}^2$.  The computation of $\tilde p$ and  $\tilde a$ is the second key computer result.\\

Since $H^5( S^2 \times S^2 \times I ; \R / \Z) = 0$, we can find an 
$\R / \Z$-cochain $\tilde w$ on $\widetilde T$ with
$$d\tilde w = (1/2)\tilde{p} \cup_1 \tilde{p} + (1/2)\tilde{p} \cup_2 \tilde{a}^2 + (1/4)\widetilde{A}(\widetilde{A} \cup_1 \widetilde{A}) + (1/2)x(\tilde{a}).$$ This is the basic equation, hence $(\tilde w, \tilde p, \tilde a) \in G^4(S^2 \times S^2 \times I)$. We can then add a cocycle to $\tilde w$ and assume that  $\tilde w  $ extends $(1/8)A_1A_2 + (1/2)y_4(a_1, a_2)$ on $T \times \{0\}$.  Thus $(\tilde w, \tilde p, \tilde a)$ represents the product class we are studying, and restricts to some appropriate $(w,0,c)$ on $T' \times \{1\}$.\\  

By Stokes' Theorem, the integral of $\tilde{w}\vert_{S^2 \times S^2 \times \{1\}}$ minus the integral of $\tilde{w}\vert_{S^2 \times S^2 \times \{0\}}$ is equal to the integral of
$$(1/2)\tilde{p} \cup_1 \tilde{p} + (1/2)\tilde{p} \cup_2 \tilde{a}^2 + (1/4)\widetilde{A}(\widetilde{A} \cup_1 \widetilde{A}) + (1/2)x(\tilde{a})$$ over $S^2 \times S^2 \times I$. Computer computation shows that this latter integral is 0. This is the third key computer  result.\\

Thus the integrals of the restrictions of $\tilde{w}$ to the 0-end and the 1-end of $S^2 \times S^2 \times I$ are
equal. The restriction of $(\tilde w, \tilde p, \tilde a)$  to $(T \times  \{0\})$ agrees identically with the product $(0,0,a_1)(0,0,a_2)$.   Therefore the integrals of $\tilde{w}$ restricted to both ends of $S^2 \times S^2 \times I$ must be $1/8$. Since the product $(0,0,a_1)(0,0, a_2)$ evaluates zero on the identity map
and since $(w, 0, c) \equiv (0,0,a_1)(0,0,a_2) \in G^4(S^2 \times S^2)$, it follows that $(0,0,c)$ evaluates
$-1/8$ on the identity map.\\

This completes the proof of the claim, and hence establishes the validity of our general formulas for the evaluation $G^4(X) \to \rm{Hom}(\Omega^{spin}_4(X), \R / \Z)$ described in \S3. Note that if we had made the alternate choice of $y_4(a, b)$ in \S5.1, by adding $ab$, then the class $w$ would have integrated to $5/8$ over $S^2 \times S^2 \times \{1\}$ and the class $(0,0, c) = (0, 0, g^*s\alpha)$ would necessarily have evaluated as  $-5/8$.  The resulting product and evaluation on $G^4(X)$ would then have been defined so as to be suspension compatble with the alternate evaluation on $G^3(X)$ discussed at the end of  \S2.2.

\section {Appendices}

\subsection {Cochain Formulas For $y_4(a, b)$ and $z(r)$}

We state here explicit formulas for the term $y_4(a, b)$ that occurs in the product formula for triples, as discussed in \S5.1, \S6.1, and \S6.2, and for the term $z(r)$ that occurs in the third relation between triples, as discussed in \S5.2 and \S7.1.  It is to be emphasized that these formulas bring to the forefront how complicated our explicit combinatorial description of the group $G^4(X)$ really is.  Of course both $z(r)$ and $y_4(a, b)$ are only defined up to natural coboundaries, but we doubt that  there are much shorter formulas.  We worked out a formula for $z(r)$ by hand. We used a computer to find $y_4(a, b)$ by solving the equation $dy_4(a, b) = \Delta x(a, b) + \delta x(a,b)$, as explained in \S5.1.\\

First, here is a formula for the term $z_M(r)$ discussed in \S5.2.  We write out the evaluation of $z_M(r)$ on a 4-simplex $(01234)$.
$$z_M(r) = r(01)r(01)r(12)r(34) + r(01)r(01)r(13)r(34) + r(01)r(12)r(23)r(24)$$  $$+ r(01)r(12)r(24)r(24) + r(01)r(12)r(24)r(34) + r(01)r(12)r(14)r(23)$$ $$+ r(01)r(12)r(14)r(24) + r(01)r(12)r(14)r(34) + r(02)r(23)r(23)r(34).$$
Some of the quartic terms in this sum are actually cubic terms as functions because in $\ZZ/2$ arithmetic $e^2 = e$.\\

\noindent Then from \S5.2  $$z(r) = z_M(r) + r(dr \cup_1dr).$$  The second term $ r(dr \cup_1dr)(01234) = r(01)(dr \cup_1 dr)(1234)$ becomes a sum of  cubic and quadratic terms in evaluations of $r$ on edges of a 4-simplex.
$$r(01)(dr \cup_1 dr)(1234) = r(01)[dr(124)dr(234) + dr(134)dr(123)] $$ $$ = r(01)[(r(12)+r(14)+r(24))(r(23)+r(24)+r(34))] $$ $$ +r(01)[ (r(13)+r(14)+r(34))(r(12)+r(13)+ r(23))].$$

Next, here comes our explicit formula for $y_4(a, b)$. Since $a$ and $b$ are $\Z/2$ cocycles, if $0 \leq i < j < k \leq 3$ one can write an evaluation $a(ijk) = a(ij4) + a(ik4) + a(jk4)$, and similarly for $b$.  This explains why every evaluation term below involves vertex 4.  In fact, in the process of finding $y_4(a, b)$, with specified $dy_4(a, b)$, it was necessary to work with a basis of the space of  $\Z/2$ 2-cocycles  on a 4-simplex, and then polynomial functions in two variables on that space. An element of the  six dimensional space of 2-cocycles on a 4-simplex is determined by naming arbitrary values $a(ij4)$, for $0 \leq i < j \leq 3$.\\

We will write out the evaluation of $y_4(a, b)$ on a 4-simplex as a sum of 174 terms on the next two pages. It may look rather silly, but we used this exact formula, albeit internally in a computer program, for an important computation in \S9.3.  The class $y_4(a, b)$  was characterized implicitly up to coboundaries in \S5.1.  But we wanted to record an explicit product formula for our group $G^4(X)$.  Looking forward, we believe that operad techniques will place our face evaluation formulas for $x(a), y_4(a, b)$, and $z(r)$ in a more conceptual context, very much like the cup$_i$ products.  These products also have formulas in terms of sums of products of face evaluations, but their theoretical characterizations and formal properties are what are considered most important.
$$y_4(a, b)(01234) = (\rm{quadratic\ terms}) + (\rm{cubic\ terms}) + (\rm{quartic\ terms})$$
where the quadratic $ab$ terms are\\

\noindent$ a(014)b(234) + a(024)b(234)+a(124)b(234) +a(014)b(014) +a(024)b(024)+a(034)b(014)+a(124)b(124)+a(134)b(014).$\\

\noindent The cubic terms of form $aab$ are\\

\noindent$ a(014)a(024)b(234)+a(014)a(034)b(024)+a(014)a(034)b(034)+a(014)a(034)b(124)+ a(014)a(034)b(234)+ a(014)a(124)b(134)+ a(014)a(124)b(234)+ a(014)a(134)b(124)+ a(014)a(134)b(134)+ a(014)a(134)b(234)+ a(014)a(234)b(134)+ a(014)a(234)b(234)+ a(024)a(034)b(034)+ a(024)a(034)b(124)+ a(024)a(134)b(124)+ a(024)a(134)b(234)+ a(024)a(234)b(234) +  a(034)a(124)b(014)+ a(034)a(124)b(024)+ a(034)a(124)b(124)+ a(034)a(124)b(134)+ a(034)a(124)b(234)+ a(034)a(134)b(014)+ a(034)a(134)b(034)+ a(034)a(234)b(014)+ a(034)a(234)b(034)+ a(034)a(234)b(134)+ a(124)a(134)b(124) $\\

\noindent The cubic terms of form $abb$ are\\

\noindent$a(014)b(014)b(234)+a(014)b(024)b(034)+a(014)b(024)b(234)+a(014)b(034)b(124)+a(014)b(034)b(134)+a(014)b(124)b(134)+a(014)b(124)b(234)+a(014)b(134)b(234)+ a(024)b(014)b(134)+ a(024)b(024)b(034)+ a(024)b(024)b(234)+ a(024)b(034)b(124)+ a(024)b(034)b(134)+ a(024)b(124)b(134)+ a(024)b(124)b(234)+ a(024)b(134)b(234) +a(034)b(014)b(024)+a(034)b(014)b(034)+a(034)b(014)b(134)+a(034)b(024)b(034)+a(034)b(024)b(124)+a(034)b(034)b(124) + a(124)b(014)b(134)+ a(124)b(014)b(234)+ a(124)b(034)b(124)+ a(124)b(034)b(134)+ a(124)b(034)b(234) +a(134)b(014)b(024)+a(134)b(014)b(034)+a(134)b(014)b(134)+a(134)b(014)b(234)+a(134)b(024)b(034)+a(134)b(024)b(124)+a(134)b(024)b(234)+a(134)b(034)b(124)+ a(234)b(014)b(234)$\\

\noindent The quartic terms of form $aaab$ are\\

\noindent$a(014)a(024)a(234 )b(014 )+a(014)a(024)a(234 )b(024 )+a(014)a(024)a(234 )b(234 )+a(014)a(034)a(234 )b(024 )+a(014)a(034)a(234 )b(034)+a(014)a(034)a(234 )b(234)+a(014)a(124)a(234 )b(124)+a(024)a(034)a(234 )b(024)+a(024)a(034)a(234 )b(034 )+a(024)a(034)a(234 )b(234 )+a(024)a(124)a(234 )b(014)+a(024)a(124)a(234 )b(024)+a(024)a(124)a(234 )b(124)+a(024)a(124)a(234 )b(234)+a(034)a(124)a(234 )b(024)+a(034)a(124)a(234 )b(034)+
a(034)a(124)a(234 )b(234)$\\

\noindent The quartic terms of form $aabb$ are\\

\noindent$a(014)a(024)b(014)b(134)+a(014)a(024)b(024)b(234)+a(014)a(024)b(034)b(124)+a(014)a(024)b(034)b(134)+a(014)a(024)b(034)b(234)+a(014)a(024)b(124)b(134)+a(014)a(024)b(124)b(234)+a(014)a(024)b(134)b(234)+a(014)a(034)b(014)b(134)+a(014)a(034)b(024)b(234)+a(014)a(034)b(034)b(234)+a(014)a(124)b(014)b(134)+a(014)a(124)b(014)b(234)+a(014)a(124)b(034)b(124)+a(014)a(124)b(034)b(134)+a(014)a(124)b(034)b(234)+a(014)a(124)b(124)b(134)+a(014)a(124)b(124)b(234)+a(014)a(124)b(134)b(234)+a(014)a(134)b(014)b(124)+a(014)a(134)b(014)b(234)+a(014)a(134)b(034)b(124)+a(014)a(134)b(034)b(134)+a(014)a(134)b(034)b(234)+a(014)a(234)b(014)b(234)+a(014)a(234)b(024)b(234)+a(014)a(234)b(124)b(234)+a(024)a(034)b(014)b(024)+a(024)a(034)b(014)b(034)+a(024)a(034)b(014)b(124)+a(024)a(034)b(014)b(134)+a(024)a(034)b(024)b(034)+a(024)a(034)b(024)b(124)+a(024)a(034)b(024)b(234)+a(024)a(034)b(034)b(134)+a(024)a(034)b(124)b(134)+a(024)a(034)b(124)b(234)+a(024)a(034)b(134)b(234)+a(024)a(124)b(014)b(234)+a(024)a(124)b(024)b(234)+a(024)a(124)b(124)b(234)+a(024)a(134)b(014)b(124)+a(024)a(134)b(014)b(134)+a(024)a(134)b(014)b(234)+a(024)a(134)b(034)b(124)+a(024)a(134)b(034)b(134)+a(024)a(134)b(034)b(234)+a(024)a(134)b(124)b(134)+a(024)a(134)b(134)b(234)+a(024)a(234)b(014)b(234)+a(024)a(234)b(124)b(234)+a(034)a(124)b(014)b(124)+a(034)a(124)b(014)b(134)+a(034)a(124)b(014)b(234)+a(034)a(124)b(024)b(234)+a(034)a(124)b(034)b(124)+a(034)a(124)b(034)b(134)+a(034)a(124)b(124)b(134)+a(034)a(124)b(134)b(234)+a(034)a(134)b(014)b(124)+a(034)a(134)b(014)b(134)+a(034)a(134)b(014)b(234)+a(034)a(134)b(034)b(124)+a(034)a(134)b(034)b(134)+a(034)a(134)b(034)b(234)+a(034)a(134)b(124)b(134)+a(034)a(134)b(134)b(234)+a(124)a(134)b(014)b(124)+a(124)a(134)b(014)b(134)+a(124)a(134)b(014)b(234)+a(124)a(134)b(034)b(124)+a(124)a(134)b(034)b(134)+a(124)a(134)b(034)b(234)+a(124)a(134)b(124)b(134)+a(124)a(134)b(124)b(234)+a(124)a(134)b(134)b(234)+a(124)a(234)b(014)b(234)+a(124)a(234)b(024)b(234)$\\

\noindent The quartic terms of form $abbb$ are\\

\noindent$a(014)b(014)b(124)b(234)+a(014)b(024)b(124)b(234)+a(034)b(014)b(024)b(234)+a(034)b(014)b(034)b(234)+a(034)b(024)b(034)b(234)+a(034)b(024)b(124)b(234)+a(034)b(034)b(124)b(234)$

\subsection {The Group Extensions for the Filtration Quotients}

We want to use our cochain and cocycle constructions to give information about the extensions that arise from the $E_\infty$ page of the Atiyah-Hirzebruch spectral sequence for $4D$ Spin bordism, or more precisely, for its Pontrjagin dual.  These extensions were mentioned in \S7.5.  The ultimate goal is, for a finite complex $X$, how might one actually express $G^4(X)$ as a  direct sum of finite cyclic groups and copies of $\R / \Z$?  If we were ambitious enough, we could probably do this.  But we will only state some partial results about the extensions below.
$$(I) \hspace{.5in}0 \to QH^4(X;\Ar/\Zee) \to G^1(X) \to SH^3(X;\Zee/2\Zee) \to 0 \hspace{.5in}$$
$$ (II)\hspace{.2in} 0 \to SH^3(X;\Zee/2\Zee)  \to G(X)/ G^2(X) \to SSH^2(X;\Zee/2\Zee) \to 0 \hspace{.2in}$$
$$ (III) \hspace{.7in} 0 \to G^1(X) \to G(X) \to SSH^2(X;\Zee/2\Zee) \to 0 \hspace{,7in}$$
$$ (IV)\hspace{.5in}  0 \to QH^4(X;\Ar/\Zee) \to G(X) \to G(X) / G^2(X) \to 0 \hspace{.6in}$$
Isomorphism classes of extensions of abelian groups can be characterized in terms of certain collections of auxillary characteristic functions. Specifically, in the simplest case, if $$0 \to K \to E \to V \to 0$$ is an extension of abelian groups and if $V$ is a $\Z/2$ vector space, then the isomorphism class of the extension is determined by the homomorphism $e \colon V \to K / 2K$ given by $e(v) = \hat{v}^2 \in K / 2K$, where $\hat{v} \in E$ lifts $v \in V$.  This fact can be applied to the first three sequences above.  If $V$ is a $\Z/4$ module the extension is determined up to isomorphism by the two homomorphisms $e(v) = \hat{v}^2 \in K / 2K$ if $2v = 0$ and $e'(v) = \hat{v}^4 \in K / 4K$ for all $v$.  This fact can be applied to the last sequence above, since from the second sequence $G(X) / G^2(X)$ is a $\Z/4$ module.\\

In fact, the second extension $(II)$ above is the easiest to understand.  One sees that $G(X)/ G^2(X)$ can be described as certain equivalence classes of pairs $(p, a)$.  Specifically, there should exist a lift $(w,p,a) \in G(X)$, so $\bar{a} \in SSH^2(X; \Z/2)$, and there are restrictions on $p$ beyond $dp = a^2$.  The product formula and relations for the pairs $(p, a) \in G(X) / G^2(X)$ are easily found from the product formula and relations in $G(X)$.  The characteristic homomorphism $e(\bar{a})$ for the second extension $(II)$ above is then computed from $ (p,a)^2 = (Sq^1a, 0)$, so $e(\bar a) =Sq^1\bar{a} \in SH^3(X; \Z/2)$.  Then $G(X) / G^2(X)$ is a direct sum of $\Z/2$'s and $\Z/4$'s, with the number of $\Z/4$ summands equal to the rank of the map $e(\bar{a}) = Sq^1\bar{a}$, $ \bar{a} \in SSH^2(X; \Z/2)$.\\

For the first sequence $(I)$,  it is helpful to use Bocksteins  of coefficient sequences  to identify $$QH^4(X;\Ar/\Zee)/ 2QH^4(X;\Ar/\Zee) \simeq (T / 2T) \simeq Image(T)  \subset H^5(X; \Z/2),$$ where $T \subset H^5(X; \Z)$ is the torsion subgroup.\footnote{For $\Z/2$ 2-cocycles $q$, the fact that elements $(1/2)q^2 = 0 \in QH^4(X;\Ar/\Zee)$ becomes irrelevant, because $(1/2)q^2 = 2(1/4)\mathcal{P}(q) \in H^4(X; \R/\Z)$.} Lift $\bar{p} \in  SH^3(X;\Zee/2\Zee)$ to $(w,p,0) \in G^1(X)$.  Then $$(w,p,0)^2 = (2w + (1/2)p \cup_2p,\ 0,\ 0) \equiv (2w, 0, 0),$$ where the last equivalence follows from the fact that $(1/2)Sq^1p$ is a coboundary in $C^4(X; \R / \Z)$.  Since $dw = (1/2)Sq^2p$, under the Bockstein identification above we have $e(\bar{p}) = Sq^2\bar{p}\in Image(T)\subset H^5(X; \Z/2)$. For spaces $X$ with finitely generated homology it still requires a little manipulation with  Bockstein sequences to express the extension $G^1(X)$ as a sum of $\R / \Z$'s and cyclic abelian groups.\\

In order to describe the characteristic homomorphism for the third sequence $(III)$, it is necessary to give some description of the $\Z/2$ vector space $G^1(X) / 2G^1(X)$.  One can see that there is a natural short exact sequence of $Z/2$ vector spaces
$$0 \to Image(T) / Sq^2SH^3(X; \Z/2) \to G^1(X) / 2G^1(X) \to SH^3(X; \Z/2) \to 0,$$ where the group on the left is a quotient of the subgroup  $Image(T) \subset H^5(X; \Z/2)$.   Of course this sequence of $\Z/2$ vector spaces splits, but it does not split  naturally.\\

Returning to the characteristic homomorphism for the extension $(III)$, lift $\bar{a} \in SSH^2(X; \Z/2)$ to $(w,p,a) \in G(X)$.  Then after some simplification of the product formula in $G(X)$ one has 
$$(w,p,a)^2 \simeq (2w -(1/8)(A^2 + 2A \cup_1(A \cup_1A)),\  Sq^1a, \ 0) \in G^1(X).$$
This element is interpreted as  $ e(\bar{a}) \in G^1(X) / 2G^1(X)$, but it is not so clear how to relate this to  the exact sequence of $\Z/2$ vector spaces above, other than the obvious remark that $Sq^1\bar{a}$ is the projection of the element $e(\bar{a})$ to $SH^3(X; \Z/2)$.\\

It is somewhat more complicated to discuss the characteristic homomorphisms $e, e'$ for the fourth extension $(IV)$ above, so we will skip this.\\

\subsection{A Subgroup of the Dual of $n$-Dimensional Spin Bordism}

Consider the two stage Postnikov tower
$$ \tilde{E} = K(\Z/2, n-1) \ltimes_{(1/2)Sq^2p} K(\R / \Z, n).$$ 
Here, $p$ is the fundamental $\Z/2$ cocycle of degree $n-1$ and $Sq^2p = p \cup_{n-3} p$ is the standard cocycle representative of the cohomology operation $Sq^2$ in this dimension. The space $\tilde{E}$ represents the Pontrjagin dual $\tilde{G}(X)$ of a quotient of reduced $n$-dimensional Spin bordism.  Specifically, from \S1.5, elements of $\tilde{G}(X) = [X, \tilde E]$ are represented by pairs $(w, p)$ with $dp = 0$ and $dw = (1/2)p \cup_{n-3} p$.  The product is given by
$$(w,p)(v,q) = (w+v + (1/2)p \cup_{n-2} q,\  p+q).$$
The null-homotopic pairs are $(df + (1/2) Sq^2c,\ dc)$, where $Sq^2c = c \cup_{n-3}dc + c  \cup_{n-4} c$.  We assert the following result, which is similar to discussions in the physics paper [4]. \\

\noindent CLAIM 1: There is a natural isomorphism $$\tilde{G}(X) \to {\rm Hom}(\widetilde{\Omega}_n^{spin}(X) / Image( \widetilde{\Omega}_n^{spin}(X^{(n-2)})),\ \R /\Z)$$ where $X^{(n-2)}$ denotes the $n-2$ skeleton of $X$.\\

To prove Claim 1, we need to evaluate a pair $(w, p)$ on a reduced Spin bordism representative $f\colon M^n \to X$. We first find a cochain $c$ on $M$ and a simplicial map $u\colon M \to S^{n-1}$ with $u^*(z) = f^*p + dc $, where $z$ is a standard cocycle representing the generator of the cohomology of the sphere.  Let $Z \subset M$ be a framed, hence Spin, 1-submanifold obtained as the transverse inverse image under $u$ of a point in $S^{n-1}$.  So $[Z] \in \Omega^{spin}_1(pt) = \Z/2$ is just another name for the bordism class $[u] \in \widetilde{\Omega}_n^{spin}(S^{n-1}) = \Z/2$. Then we set $$\langle (w,p), (M \xrightarrow{f} X) \rangle = $$ $$ (1/2)[Z]  + (1/2)\int_{[M]} Sq^2c + (1/2)\int_{[M]} f^*p \cup_2 dc + \int_M f^*w \in \R / \Z.$$  This construction has already appeared in \S3.4 in the case $n = 4$, and also appears in our paper ([1], \S6) in the case $n = 3$.\\

We can also use Kapustin's method described in \S3.4 to define  evaluations when $n \geq 3$. We can assume $X = M$ and $f = Id$. The reduced Spin bordism in dimension $n$ of products of $K(\Z / 2, n-1)$'s vanishes.   Choose a Spin manifold  $W^{n+1}$ with $\partial W = M$, so that all cocycles $p \in Z^{n-1}(M; \Z/2)$ lift to cocycles $\tilde p \in Z^{n-1}(W ; \Z/2)$.  Then define
$$\langle (w,p), (M \xrightarrow{Id} M) \rangle =   (1/2)\int_W\ Sq^2\tilde p + \int_M\ w \in \R / \Z.$$

It is not at all obvious that the two evaluations coincide. Before proving this, we establish some properties of the Kapustin evaluation. From the fact that cohomologically $Sq^2$ vanishes into the top dimension for closed Spin manifolds or Spin manifold pairs $(W, \partial W)$, it is not hard to see that the function $Q(p) = \int_W\ Sq^2 \tilde p \in /Z/2$ is well-defined, that is, independent of the choice of $W$ and $\tilde p$.   Next, we prove that the Kapustin evaluation vanishes on pairs $(df + (1/2)Sq^2c, dc)$ that represent zero in $\tilde{G}(M)$. The $df$ causes no trouble since $\int_{[M]}\ df = 0$. Lift cochain $c$ to a cochain $\tilde c$ on $W$.  Then $$Q(dc) = \int_{[W]}\ Sq^2(d\tilde c) = \int_{[W]}\ dSq^2(\tilde c) = \int_{[M]}\ Sq^2(c),$$ which is what we want.  Also,  the function $Q$ is  quadratic over the pairing $\langle p \cup_{n-2} q, [M] \rangle$. This means 
$$Q(p+q) = Q(p) + Q(q) + \int_{[M]}\ p\cup_{n-2} q.$$  The quadratic statement follows from the basic property of cup$_i$ products  $$(\tilde p+ \tilde q) \cup_{n-3} (\tilde p+ \tilde q) = \tilde p \cup_{n-3} \tilde p + \tilde q \cup_{n-3}\tilde q + d(\tilde p \cup_{n-2} \tilde q),$$ along with Stokes Theorem. The quadratic property of $Q$ easily implies  that the Kapustin evaluation defines a group homomorphism   $$\tilde{G}(X) \to \rm{Hom}(\widetilde{\Omega}_n^{spin}(X),\ \R /\Z).$$  The image obviously vanishes on $Image( \widetilde{\Omega}_n^{spin}(X^{(n-2)})) \to \widetilde{\Omega}_n^{spin}(X))$.\\

In ([1], \S6.3) we gave an argument that the two evaluations coincide for $n = 3$.  Here is an easier argument for $n > 3$.  A little computation using the quadratic property for $Q$ and the formula  $Q(dc) = \int_{[M]} Sq^2c$ shows that it suffices to prove\footnote{Again, $[Z]$ denotes the Spin bordism class of a framed 1-manifold $Z = u^{-1}(pt) \subset M$. We interpret both $Q(p) = \int_{W]}Sq^2\tilde p$ and $[Z]$ as elements of $\Z/2$.} $Q(p) = [Z]$, if $p = u^*(z)$ for a map $u \colon M \to S^{n-1}$. The reduced $n$-dimensional Spin bordism of $\Sigma^{n-3}\mathbb{C}P(2)$ vanishes, so we can extend $u$ to a map $\tilde u \colon (W, M) \to (\Sigma^{n-3}\mathbb{C}P(2), S^{n-1})$, where $W$ is a Spin manifold with $\partial W = M$.  Then $\tilde u^*(\tilde z) = \tilde p$, where $\tilde z$ is a cocycle representing the generator of $H^{n-1}(\Sigma^{n-3} \mathbb{C}P^2; \Z/2)$.  We have $ Sq^2 \tilde z \not = 0$, and this cocycle $Sq^2 \tilde z = \tilde z \cup_{n-3} \tilde z$ is a relative cocycle on the pair $(\Sigma^{n-3}\mathbb{C}P(2), S^{n-1})$.  This is the key. On the one hand, in the top dimension $$\int_{[W]}\ \tilde{u}^*Sq^2(\tilde z) = \int_{[W]}\ Sq^2 \tilde p = Q(p).$$  On the other hand, this integral is exactly the obstruction to deforming $\tilde u \colon (W, M) \to (\Sigma^{n-3}\mathbb{C}P(2), S^{n-1})$ rel $M$ to a map $W \to S^{n-1}$, and this obstruction is the Spin bordism class of $u\colon M \to S^{n-1}$, also known as $[Z]$.  Thus the two evaluations coincide. \\

From the relations on representative pairs  $(w, p) \in \widetilde{G}(X)$,  we see that there is a short exact sequence\footnote{The characteristic homomorphism for this sequence  is given by $e(\bar{p}) = Sq^2(\bar{p}) \in T/(2T+Image(\beta Sq^2)) \subset H^{n+1}(X; \Z/2)/Image(Sq^3)$, where $T \subset H^{n+1}(X; \Z)$ is the torsion subgroup and $T / (2T + Image(\beta Sq^2)) \simeq  QH^n(X; \R/\Z) / 2QH^n(X; \R/\Z) $.  This claim basically follows in a manner similar to the   discussion of the extension $(I)$ in Appendix 10.2. }
$$0 \to QH^n(X; \R /\Z) \to \tilde{G}(X) \to SH^{n-1}(X; \Z/2) \to 0,$$
where $$QH^n(X; \R /\Z) = H^n(X; \R / \Z)\ /\ Image(H^{n-2}(X; \Z/2) \xrightarrow{(1/2)Sq^2} H^n(X; \R / \Z))$$ and  $$SH^{n-1}(X; \Z/2) = Kernel( H^{n-1}(X; \Z/2) \xrightarrow{(1/2)Sq^2} H^{n+1}(X ; \R / \Z)).$$
This exact sequence can also be seen by mapping $X$ into the sequence of fibrations
$$K(\Z,2, n-2) \to K(\R / \Z, n) \to \tilde E \to K(\Z/2, n-1) \to K(\R/\Z, n+1),$$ where the first and last map are both $(1/2)Sq^2$.
 Now the filtration argument given in  ([1], \S6.4) for $n = 3$ can be extended to prove that the evaluation we have defined for $n >3$ does define an isomorphism$$\tilde{G}(X) \simeq {\rm Hom}(\widetilde{\Omega}_n^{spin}(X) / Image( \widetilde{\Omega}_n^{spin}(X^{(n-2)})),\ \R /\Z).$$ The filtration argument compares the short exact sequence above involving $\tilde{G}(X)$ to the dual of the short exact sequence coming from the three quotients associated to the filtration
 $$ Image(\widetilde \Omega^{spin}_n(X^{n-2})) \subset  Image(\widetilde \Omega^{spin}_n(X^{n-1})) \subset \widetilde \Omega^{spin}_n(X).$$
For this, as in [1],  instead of the Kapustin evaluation it is easier to use the evaluation $$\langle (w,p), (M \xrightarrow{f} X) \rangle =   (1/2)[Z] + \int_M\ f^*w \in \R / \Z,$$ in the case $f^*p = u^*(z),\ u\colon M \to S^{n-1}$.  The point is, the coefficient group $\Z/2$ in $SH^{n-1}(X; \Z/2)$ really means the Pontrjagin dual of $\widetilde{\Omega}^{spin}_1(pt) = \widetilde{\Omega}^{spin}_n(S^{n-1})$.  Then a map $M \to X^{n-1} / X^{n-2}$ is a map from $M$ to a wedge of $(n-1)$-spheres.  In any case, the filtration argument completes the proof of Claim 1.\\

Now consider a Spin manifold $M^n$.  We have seen that the Kapustin evaluation above yields a canonical quadratic function $Q\colon Z^{n-1}(M; \Z/2) \to \Z/2$, satisfying
$$Q(p+q) = Q(p) + Q(q) + \int_{[M]}\  p \cup_{n-2} q $$ and $$Q(dc) = \int_{[M]}\ Sq^2c\ =  \int_{[M]}\  c \cup_{n-3} dc + c \cup_{n-4} c.$$
It is easy to see that any other such quadratic function must differ from $Q$ by a linear function vanishing on coboundaries, hence must be of form  $Q_a(p) = Q(p) + \langle ap, [M] \rangle$, for some $a \in H^1(M; \Z/2)$.  Since equivalence classes of Spin structures on $M$ are also a torsor of $H^1(M; \Z/2)$, we obtain the following quite interesting result, which can more or less be found in the physics paper [4].\\

\noindent CLAIM 2:   Let $M^n$ be a closed oriented n-manifold with a simplicial structure.  Equivalence classes of Spin structures on $M^n$ are in canonical bijective correspondence with functions $$Q\colon Z^{n-1}(M; \Z/2) \to  \Z / 2$$ that satisfy $$Q(p+q) = Q(p) + Q(q) + \int_{[M]}\  p \cup_{n-2} q $$ and $$Q(dc) = \int_{[M]}\ Sq^2c\ =  \int_{[M]}\  c \cup_{n-3} dc + c \cup_{n-4} c.$$
Such quadratic functions $Q$  exist only if the Wu class $v_2(M) = 0$.   If  $Q$ is the canonical quadratic function associated to one Spin structure on $M$, and if  $M_a$ denotes the new Spin structure on the same oriented manifold, obtained by acting by  $a \in H^1(M; \Z/2)$, then the Spin-1 manifold evaluation  of $Q$ can be used to show that the new canonical quadratic function on $M$ is exactly $Q_a$.\\

In a subsequent paper, [2], we will study these quadratic functions  in greater detail.  In particular, we will extend the construction to Spin manifolds with boundary, and we will show how quadratic functions on manifolds directly  induce quadratic functions on boundaries, and also on codimension 0 submanifolds.  These constructions make use of the cochain suspension operation $s$ that we introduced in \S8 of this paper.

\section* {\bf References}

1. Greg Brumfiel and John Morgan, The Pontrjagin Dual of 3-Dimensional Spin Bordism, arXiv.org$>$math.AT$>$arXiv:1612.02860v1.\\

\noindent2.  Greg Brumfiel and John Morgan, Quadratic Functions of Cocycles and Spin Structures, to appear.\\

\noindent3. R. R. Bruner and J. P. C. Greenlees, The connective real K-theory of finite groups,
Math. Surveys and Monographs 169, 2010.\\

\noindent4. D. Gaiotto and A. Kapustin,  Spin TQFTs and fermonic phases of matter,
arXiv.org$>$cond-mat$>$arXiv:1505.05856v2.\\ 

\noindent5. Z-C. Gu and X-G. Wen, Symmetry-protected topological orders for interacting fermions: Fermionic topological nonlinear $\sigma$-models and a special group supercohomology theory, Phys. Rev. B 90, 115141 (2014) [arXiv:1201.2648 [cond-mat]\\

\noindent 6 . A.  Kapustin,  Symmetry  Protected  Topological  Phases, Anomalies,  and Cobordisms:  Beyond Group Cohomology,
arXiv:1403.1467 [cond-mat.str-el].

\end{document}